\documentclass[a4paper,twoside,11pt]{article}
\usepackage{amsfonts}
\usepackage{amssymb}\usepackage{amsmath}
\usepackage{graphicx}
\def\be{\begin{equation}}
\def\ee{\end{equation}}
\def\ba{\begin{array}}
\def\ea{\end{array}}
\def\bc{\begin{center}}
\def\ec{\end{center}}

\def\ZZ{\rm {{\rm Z}\kern-.48em{\rm Z}}}
\def\RR{\rm \hbox{I\kern-.2em\hbox{R}}}
\def\CC{\rm \hbox{C\kern -.5em{\raise .32ex
\hbox{$\scriptscriptstyle |$}}\kern - .22em{\raise .6ex
\hbox{$\scriptscriptstyle |$}}\kern
.4em}}
\def\CC{ {\mathbf C} }
\def\RR{ {\mathbf R} }
\def\ZZ{ {\mathbf Z} }

\newtheorem{Tm}{Theorem}[section]
\newtheorem{Df}[Tm]{Definition}

\newtheorem{Lm}[Tm]{Lemma}

 \numberwithin{equation}{section}
\newcommand{\double}{\baselineskip 1.24 \baselineskip}

\title{ $(p,q)$-frames in shift-invariant subspaces of mixed Lebesgue spaces $L^{p,q}(\mathbf{R}\times \mathbf{R}^{d})$}

\author{{\small  Yingchun Jiang and Jiao Li }\\
{\small School of Mathematics and Computational Science,}\\{\small
Guilin University of Electronic Technology, Guilin, 541004, P. R. China}}

\begin{document}
\date{}
\maketitle \double
\noindent\textbf{Abstract:}\  In this paper, we mainly discuss the $(p,q)$-frames in shift-invariant subspaces
 \begin{equation*}
 V_{p,q}(\Phi)=\left\{\sum\limits_{i=1}^{r}\sum\limits_{j_{1}\in \mathbf{Z}}\sum\limits_{j_{2}\in \mathbf{Z}^{d}}d_{i}(j_{1},j_{2})\phi_{i}(\cdot-j_{1},\cdot-j_{2}):\Big(d_{i}(j_{1},j_{2})\Big)_{(j_{1},j_{2})\in \mathbf{Z}\times\mathbf{Z}^{d}}\in \ell^{p,q}(\mathbf{Z}\times\mathbf{Z}^d)\right\}
 \end{equation*}
 of mixed Lebesgue spaces $L^{p,q}(\mathbf{R}\times \mathbf{R}^{d})$. Some equivalent conditions for $\{\phi_{i}(\cdot-j_{1},\cdot-j_{2}):(j_{1},j_{2})\in\mathbf{Z}\times\mathbf{Z}^d,1\leq i\leq r\}$ to constitute a $(p,q)$-frame of $V_{p,q}(\Phi)$ are given. Moreover, the result shows that $V_{p,q}(\Phi)$ is closed under these equivalent conditions of $(p,q)$-frames for the family $\{\phi_{i}(\cdot-j_{1},\cdot-j_{2}):(j_{1},j_{2})\in\mathbf{Z}\times\mathbf{Z}^d,1\leq i\leq r\}$, although the general result is not correct.

 \noindent{\bf MR(2000) Subject Classification:}  94A20, 46E22, 62D05.

\section{Introduction }
\label{intro}
\ \ \ \ \ Frames were first introduced by Duffin and Schaeffer in the context of nonharmonic Fourier series $\cite{R1}$, and were furtherly studied in many references,  such as $\cite{A1,C2,C3,H1,K2,R2,R3,W2}$. In recent years, frames have been generally applied  to wavelet theory, time frequency analysis and sampling theory $\cite{C1,D2,D1,A3}$.
\par Some signals are time-varying in practice, which means that the signals live in time-space domains at the same time. Mixed Lebesgue space is a suitable tool for modeling and measuring time-space signals, due to the separate integrability for different variables. Mixed Lebesgue
spaces were first described in detail by Benedek and Panzone $\cite{A4}$, and were furtherly studied from the views of classical harmonic analysis and operator theory, see $\cite{A5,G1,D3,J1,A6,L4,J2}$ and the references therein.
\par For $1\leq p,q\leq \infty$, the mixed Lebesgue space $L^{p,q}(\mathbf{R}\times \mathbf{R}^{d})$ denotes the Banach space of all functions $f$ such that
\begin{equation}
\|f\|_{L^{p,q}(\mathbf{R}\times \mathbf{R}^{d})}=\Big\|\big\|f(x_1,x_2)\big\|_{L^q_{x_2}(\mathbf{R}^d)}\Big\|_{L^p_{x_1}(\mathbf{R})}<\infty.
\end{equation}
Similarly, $\ell^{p,q}=\ell^{p,q}(\mathbf{Z}\times \mathbf{Z}^{d})$ is the Banach space of all sequences $d=\Big(d(j_{1},j_{2})\Big)_{(j_{1},j_{2})\in \mathbf{Z}\times \mathbf{Z}^{d}}$ such that
\begin{equation}
\|d\|_{\ell^{p,q}}=\Big\|\big\|d(j_1,j_2)\big\|_{\ell^q_{j_2}(\mathbf{Z}^d)}\Big\|_{\ell^p_{j_1}(\mathbf{Z})}<\infty.
\end{equation}
Recently,  sampling for time-varying signals in bandlimited subspaces, shift-invariant subspaces and reproducing kernel subspaces of mixed Lebesgue spaces $L^{p,q}(\mathbf{R}\times \mathbf{R}^d)$ are studied in $\cite{Y1,Y2,K1,L2,L3,Q1}$. Moreover, the $L^{p,q}$-stability of the shifts of finitely many functions in mixed Lebesgue spaces $L^{p,q}(\mathbf{R}^{d+1})$ was discussed in $\cite{L1}$, which generalized the corresponding results of $L^{p}$-stability in \cite{R5, R4,R3}.

In this paper, we mainly  study the $(p,q)$-frames of the form $\big\{\phi_{i}(\cdot-j_{1},\cdot-j_{2}):(j_{1},j_{2})\in \mathbf{Z}\times\mathbf{Z}^{d},1\leq i\leq r\big\}$ in shift-invariant subspaces $V_{p,q}(\Phi)$ of $L^{p,q}(\mathbf{R}^{d+1})$, which  is  a direct generalization of the notion of $p$-frames in $\cite{A1}$.

\begin{Df}
The family $\big\{\phi_{i}(x_{1}-j_{1},x_{2}-j_{2}):\;i=1,\cdots,r,\;(j_{1},j_{2})\in \mathbf{Z}\times\mathbf{Z}^{d}\big\}$ is called a $(p,q)$-frame of $V_{p,q}(\Phi)$ if $\phi_{i}\in L^{p',q'}(\mathbf{R}\times \mathbf{R}^{d})$ and there exist positive constants $A$ and $B$ such that for any $f\in L^{p,q}(\mathbf{R}\times \mathbf{R}^d)$,
\begin{equation}
A\|f\|_{L^{p,q}}\leq\sum\limits_{i=1}^{r}\Big\|\Big\{<f,\phi_{i}(\cdot-j_{1},\cdot-j_{2})>\Big\}_{(j_{1},j_{2})\in \mathbf{Z}\times \mathbf{Z}^{d}}\Big\|_{\ell^{p,q}}\leq B\|f\|_{L^{p,q}},
\end{equation}
where $p'$ and $q'$ are the conjugate numbers of  $p$ and $q$, respectively.
\end{Df}
\par The paper is organized as follows. In section $2$, we give some lemmas which will be used for proving our main results in the next section.  The main theorem and all the proofs are gathered in section $3$.

\section{Some lemmas}

\ \ \ \ In this section, we will establish some lemmas which are useful for proving the main results. Firstly, we introduce the function spaces
\begin{equation*}
\mathcal{L}^{p,q}=\bigg\{ f:
\|f\|_{\mathcal{L}^{p,q}}=\bigg\|\sum \limits_{j_1\in \mathbf{Z}}\Big\|\sum\limits_{j_{2}\in\mathbf{Z}^{d}}|f(x_{1}+j_{1},x_{2}+j_{2})|\Big\|_{L^{q}_{x_2}([0,1]^d)}\bigg\|_{L^{p}_{x_1}([0,1])}<\infty\bigg\},
\end{equation*}
\begin{equation*}
\mathcal{W}(L^{1,1})=\bigg\{f:\|f\|_{\mathcal{W}(L^{1,1})}=\sum\limits_{j_{1}\in\mathbf{Z}}\sup\limits_{x_{1}\in [0,1]}\sum\limits_{j_{2}\in \mathbf{Z}^{d}}\sup\limits_{x_{2}\in [0,1]^{d}}|f(x_{1}+j_{1},x_{2}+j_{2})|<\infty\bigg\},
\end{equation*}
which are the generalization of $\mathcal{L}^{p}$ and $\mathcal{W}(L^{1})$ in \cite{A1} defined by
\begin{eqnarray*}
\mathcal{L}^{p}=\bigg\{f:\|f\|_{\mathcal{L}^{p}}=\Big\|\sum \limits_{j\in \mathbf{Z}^d}|f(x+j)|\Big\|_{L^{p}([0,1]^d)}<\infty\bigg\},\\
\mathcal{W}(L^{1})=\bigg\{f:\|f\|_{\mathcal{W}(L^{1})}=\sum\limits_{j\in\mathbf{Z}^{d}}\sup_{x\in[0,1]^{d}}|f(x+j)|<\infty\bigg\}.
\end{eqnarray*}
Obviously, we have $\mathcal{W}(L^{1,1})\subset\mathcal{L}^{\infty,\infty}\subset\mathcal{L}^{p,q}\subset L^{p,q}$ for $1\leq p,q\leq \infty$. For $1\leq p_{1}\leq p_{2}\leq \infty$, $\mathcal{L}^{p_{2},q}\subset\mathcal{L}^{p_{1},q}$. For $1\leq q_{1}\leq q_{2}\leq \infty$, $\mathcal{L}^{p,q_{2}}\subset\mathcal{L}^{p,q_{1}}$.

\par For any sequence $D=\Big(d(j_{1},j_{2})\Big)_{(j_{1},j_{2})\in \mathbf{Z}\times\mathbf{Z}^{d}}\in\ell^{p,q}$ and $f\in \mathcal{L}^{p,q}(\mathbf{R}\times\mathbf{R}^{d})$,  define the semi-convolution
\begin{equation*}
f\ast^{'}D=\sum_{j_{1}\in \mathbf{Z}}\sum_{j_{2}\in \mathbf{Z}^{d}}d(j_{1},j_{2})f(\cdot-j_{1},\cdot-j_{2}).
 \end{equation*}
It is easy to verify that $f\ast^{'}$ is a continuous map from $\ell^{p,q}$ to $L^{p,q}$.
\begin{Lm}\label{Lm 2.1}\cite{L1}
Let $1\leq p,q\leq \infty$. If $D\in \ell^{p,q}(\mathbf{Z}^{d+1})$ and $f\in\mathcal{L}^{p,q}(\mathbf{R}^{d+1})$, then
\begin{equation*}\label{E4}
\|f\ast^{'}D\|_{L^{p,q}}\leq \|D\|_{\ell^{p,q}}\|f\|_{\mathcal{L}^{p,q}}.
\end{equation*}
\end{Lm}
\begin{Lm}\label{Lm 2.3}\cite{A1}
Let $f\in \mathcal{L}^{p}\;(or\;\mathcal{W}(L^{1}))$, $1\leq p\leq \infty.$ Then for any $D\in \ell^{1}$,
\begin{equation*}\label{E0}
\|f\ast^{'}D\|_{\mathcal{L}^{p}}\leq \|D\|_{\ell^{1}}\|f\|_{\mathcal{L}^{p}}\;\;and\;\;\|f\ast^{'}D\|_{\mathcal{W}(L^{1})}\leq \|D\|_{\ell^{1}}\|f\|_{\mathcal{W}(L^{1})}.
\end{equation*}
\end{Lm}
\begin{Lm}\label{lm 2.4}
Let $D\in \ell^{1,1}$ and $f\in \mathcal{L}^{p,q}\;(1\leq p,q \leq\infty)$. Then
\begin{equation*}\label{E5}
\|f\ast^{'}D\|_{\mathcal{L}^{p,q}}\leq \|D\|_{\ell^{1,1}}\|f\|_{\mathcal{L}^{p,q}}.
\end{equation*}
\end{Lm}
{\bfseries Proof} Denote $d_{j_{1}}(\cdot)=d(j_{1},\cdot)$ and $f_{x_{1}}(\cdot)=f(x_{1},\cdot)$. Then
\begin{eqnarray}\label{r0}
\|f\ast'D\|_{\mathcal{L}^{p,q}}^{p}&&=\int_{[0,1]}\bigg(\sum\limits_{k_{1}\in\mathbf{Z}}\bigg(\int_{[0,1]^{d}}\bigg(\sum\limits_{k_{2}\in \mathbf{Z}^{d}}\Big|\sum\limits_{j_{1}\in\mathbf{Z}}(f_{x_{1}+k_{1}-j_{1}}\ast'd_{j_{1}})(x_{2}+k_{2})\Big|\bigg)^{q}dx_{2}\bigg)^{1/q}\bigg)^{p}dx_{1}\nonumber\\
&&=\int_{[0,1]}\bigg(\sum\limits_{k_{1}\in\mathbf{Z}}\Big\|\sum\limits_{j_{1}\in\mathbf{Z}}(f_{x_{1}+k_{1}-j_{1}}\ast'd_{j_{1}})(\cdot)\Big\|_{\mathcal{L}^{q}}\bigg)^{p}dx_{1}\nonumber\\
&&\leq \int_{[0,1]}\bigg(\sum\limits_{k_{1}\in\mathbf{Z}}\sum\limits_{j_{1}\in\mathbf{Z}}\|(f_{x_{1}+k_{1}-j_{1}}\ast'd_{j_{1}})(\cdot)\|_{\mathcal{L}^{q}}\bigg)^{p}dx_{1}\nonumber\\
&&\leq \int_{[0,1]}\bigg(\sum\limits_{k_{1}\in\mathbf{Z}}\sum\limits_{j_{1}\in\mathbf{Z}}\|d_{j_{1}}\|_{\ell^{1}}\|f_{x_{1}+k_{1}-j_{1}}\|_{\mathcal{L}^{q}}\bigg)^{p}dx_{1}\nonumber\\
&&=\int_{[0,1]}\bigg(\sum\limits_{k_{1}\in\mathbf{Z}}\sum\limits_{j_{1}\in\mathbf{Z}}\|d(j_{1},\cdot)\|_{\ell^{1}}\|f(x_{1}+k_{1}-j_{1},\cdot)\|_{\mathcal{L}^{q}}\bigg)^{p}dx_{1}.
\end{eqnarray}
Denote $c(j_{1})=\|d(j_{1},\cdot)\|_{\ell^{1}}$ and $h(x_{1})=\|f(x_{1},\cdot)\|_{\mathcal{L}^{q}}$. Then by Lemma $\ref{Lm 2.3}$
\begin{eqnarray*}
&&\int_{[0,1]}\bigg(\sum\limits_{k_{1}\in\mathbf{Z}}\sum\limits_{j_{1}\in\mathbf{Z}}\|d(j_{1},\cdot)\|_{\ell^{1}}\|f(x_{1}+k_{1}-j_{1},\cdot)\|_{\mathcal{L}^{q}}\bigg)^{p}dx_{1}\\
&&=\int_{[0,1]}\bigg(\sum\limits_{k_{1}\in\mathbf{Z}}\sum\limits_{j_{1}\in\mathbf{Z}}c(j_{1})h(x_{1}+k_{1}-j_{1})\bigg)^{p}dx_{1}\\
&&=\int_{[0,1]}\bigg(\sum\limits_{k_{1}\in\mathbf{Z}}(h\ast'c)(x_{1}+k_{1})\bigg)^{p}dx_{1}\\
&&=\|h\ast'c\|_{\mathcal{L}^{p}}^{p}\leq \|c\|_{\ell^{1}}^{p}\|h\|_{\mathcal{L}^{p}}^{p}=\|D\|_{\ell^{1,1}}^{p}\|f\|_{\mathcal{L}^{p,q}}^{p}.
\end{eqnarray*}
This together with $(\ref{r0})$ leads to $\|f\ast^{'}D\|_{\mathcal{L}^{p,q}}\leq \|D\|_{\ell^{1,1}}\|f\|_{\mathcal{L}^{p,q}}$.
\begin{Lm}\label{L2.4}
Let $D\in\ell^{1,1}$ and $f\in\mathcal{W}(L^{1,1})$. Then
\begin{equation*}
\|f\ast'D\|_{\mathcal{W}(L^{1,1})}\leq\|D\|_{\ell^{1,1}}\|f\|_{\mathcal{W}(L^{1,1})}.
\end{equation*}
\end{Lm}
{\bfseries Proof} By the similar method as in Lemma 2.3, we obtain
\begin{eqnarray*}\label{r5}
\|f\ast'D\|_{\mathcal{W}(L^{1,1})}&&=\sum\limits_{k_{1}\in\mathbf{Z}}\sup\limits_{x_{1}\in[0,1]}\sum\limits_{k_{2}\in\mathbf{Z}^{d}}\sup\limits_{x_{2}\in[0,1]^{d}}\Big|\sum\limits_{j_{1}\in\mathbf{Z}}\sum\limits_{j_{2}\in\mathbf{Z}^{d}}d_{j_{1}}(j_{2})f_{x_{1}+k_{1}-j_{1}}(x_{2}+k_{2}-j_{2})\Big|\nonumber\\
&&=\sum\limits_{k_{1}\in\mathbf{Z}}\sup\limits_{x_{1}\in[0,1]}\sum\limits_{k_{2}\in\mathbf{Z}^{d}}\sup\limits_{x_{2}\in[0,1]^{d}}\Big|\sum\limits_{j_{1}\in\mathbf{Z}}(f_{x_{1}+k_{1}-j_{1}}\ast'd_{j_{1}})(x_{2}+k_{2})\Big|\nonumber\\
&&=\sum\limits_{k_{1}\in\mathbf{Z}}\sup\limits_{x_{1}\in[0,1]}\Big\|\sum\limits_{j_{1}\in\mathbf{Z}}(f_{x_{1}+k_{1}-j_{1}}\ast'd_{j_{1}})(\cdot)\Big\|_{\mathcal{W}(L^{1})}\nonumber\\
&&\leq \sum\limits_{k_{1}\in\mathbf{Z}}\sup\limits_{x_{1}\in[0,1]}\sum\limits_{j_{1}\in\mathbf{Z}}\|d_{j_{1}}\|_{\ell^{1}}\|f_{x_{1}+k_{1}-j_{1}}\|_{\mathcal{W}(L^{1})}\nonumber\\
&&=\sum\limits_{k_{1}\in\mathbf{Z}}\sup\limits_{x_{1}\in[0,1]}\sum\limits_{j_{1}\in\mathbf{Z}}\|d(j_{1},\cdot)\|_{\ell^{1}}\|f(x_{1}+k_{1}-j_{1},\cdot)\|_{\mathcal{W}(L^{1})}\nonumber\\
&&=\sum\limits_{k_{1}\in\mathbf{Z}}\sup\limits_{x_{1}\in[0,1]}\sum\limits_{j_{1}\in\mathbf{Z}}c(j_{1})h(x_{1}+k_{1}-j_{1})\\
&&=\sum\limits_{k_{1}\in\mathbf{Z}}\sup\limits_{x_{1}\in[0,1]}(h\ast'c)(x_{1}+k_{1})\\
&&=\|h\ast'c\|_{\mathcal{W}(L^{1})}\leq \|c\|_{\ell^{1}}\|h\|_{\mathcal{W}(L^{1})}=\|D\|_{\ell^{1,1}}\|f\|_{\mathcal{W}(L^{1,1})},
\end{eqnarray*}
where $c(j_{1})=\|d(j_{1},\cdot)\|_{\ell^{1}}$ and $h(x_{1})=\|f(x_{1},\cdot)\|_{\mathcal{W}(L^{1})}$.
\begin{Lm}\label{lem 2.4}\cite{A1}
Let $1\leq p\leq\infty$. Suppose that $f\in L^{p}(\mathbf{R}^{d})$ and $g\in \mathcal{L}^{\infty}(\mathbf{R}^{d})$, then
\begin{equation*}\label{r4}
\bigg\|\bigg\{\int_{\mathbf{R}^{d}}f(x)g(x-j)dx\bigg\}_{j\in \mathbf{Z}^{d}}\bigg\|_{\ell^{p}}\leq\|f\|_{L^{p}}\|g\|_{\mathcal{L}^{\infty}}.
\end{equation*}
\end{Lm}
\begin{Lm}\label{L2.5}
Let $1\leq p,q\leq \infty.$ If $f\in L^{p,q}(\mathbf{R}\times\mathbf{R}^{d})$ and $g\in \mathcal{L}^{\infty,\infty}(\mathbf{R}\times\mathbf{R}^{d})$, then
\begin{equation*}\label{E3}
\bigg\| \left\{ \int_{\mathbf{R}}\int_{\mathbf{R}^{d}}f(x_{1},x_{2})g(x_{1}-j_{1},x_{2}-j_{2})dx_{1}dx_{2}\right\}_{(j_{1}，j_{2})\in\mathbf{Z}\times\mathbf{Z}^{d}}  \bigg\|_{\ell^{p,q}}\leq \|f\|_{L^{p,q}}\|g\|_{\mathcal{L}^{\infty,\infty}}.
\end{equation*}
\end{Lm}
{\bfseries Proof} Let $f_{x_{1}}(\cdot)=f(x_{1},\cdot)$ and $g_{x_{1}}(\cdot)=g(x_{1},\cdot)$. Then it follows from Lemma $\ref{lem 2.4}$ that
\begin{eqnarray*}
&&\bigg(\sum\limits_{j_{2}\in\mathbf{Z}^{d}}\bigg|\int_{\mathbf{R}^{d}}f_{x_{1}}(x_{2})g_{x_{1}-j_{1}}(x_{2}-j_{2})dx_{2}\bigg|^{q}\bigg)^{1/q}\\
&&=\bigg\|\bigg\{\int_{\mathbf{R}^{d}}f_{x_{1}}(x_{2})g_{x_{1}-j_{1}}(x_{2}-j_{2})dx_{2}\bigg\}_{j_{2}\in\mathbf{Z}^{d}}\bigg\|_{\ell^{q}}\\
&&\leq\|f_{x_{1}}\|_{L^{q}}\|g_{x_{1}-j_{1}}\|_{\mathcal{L}^{\infty}}.
\end{eqnarray*}
Moreover, by Minkowshi's inequality, we obtain
\begin{eqnarray*}
&&\bigg(\sum\limits_{j_{2}\in\mathbf{Z}^{d}}\bigg|\int_{\mathbf{R}}\int_{\mathbf{R}^{d}}f(x_1,x_{2})g(x_1-j_1,x_{2}-j_{2})dx_{1}dx_{2}\bigg|^{q}\bigg)^{1/q}\\
&&\leq\bigg(\sum\limits_{j_{2}\in\mathbf{Z}^{d}}\bigg(\int_{\mathbf{R}}\bigg|\int_{\mathbf{R}^{d}}f_{x_{1}}(x_{2})g_{x_{1}-j_{1}}(x_{2}-j_{2})dx_{2}\bigg|dx_{1}\bigg)^{q}\bigg)^{1/q}\\
&&\leq\int_{\mathbf{R}}\bigg(\sum\limits_{j_{2}\in\mathbf{Z}^{d}}\bigg|\int_{\mathbf{R}^{d}}f_{x_{1}}(x_{2})g_{x_{1}-j_{1}}(x_{2}-j_{2})dx_{2}\bigg|^{q}\bigg)^{1/q}dx_{1}\\
&&\leq\int_{\mathbf{R}}\|f_{x_{1}}\|_{L^{q}}\|g_{x_{1}-j_{1}}\|_{\mathcal{L}^{\infty}}dx_{1}.
\end{eqnarray*}
Finally, the desired result follows from
\begin{eqnarray*}
&&\bigg\| \left\{ \int_{\mathbf{R}}\int_{\mathbf{R}^{d}}f(x_{1},x_{2})g(x_{1}-j_{1},x_{2}-j_{2})dx_{1}dx_{2}\right\}_{(j_{1},j_{2})\in\mathbf{Z}\times\mathbf{Z}^{d}}  \bigg\|_{\ell^{p,q}}\\
&&=\bigg\| \bigg(\sum\limits_{j_{2}\in\mathbf{Z}^{d}} \bigg|\int_{\mathbf{R}}\int_{\mathbf{R}^{d}}f_{x_{1}}(x_{2})g_{x_{1}-j_{1}}(x_{2}-j_{2})dx_{1}dx_{2}\bigg|^{q}\bigg)^{1/q}\bigg\|_{\ell^{p}}\\
&&\leq\bigg\|\bigg\{\int_{\mathbf{R}}\|f_{x_{1}}\|_{L^{q}}\|g_{x_{1}-j_{1}}\|_{\mathcal{L}^{\infty}}dx_{1}\bigg\}_{j_{1}\in\mathbf{Z}}\bigg\|_{\ell^{p}}\\
&&\leq \big\|\|f_{x_{1}}\|_{L^{q}}\big\|_{L^{p}}  \big\|\|g_{x_{1}-j_{1}}\|_{\mathcal{L}^{\infty}}\big\|_{\mathcal{L}^{\infty}}\\
&&=\|f\|_{L^{p,q}}\|g\|_{\mathcal{L}^{\infty,\infty}}.
\end{eqnarray*}
The proof is completed.

 For function $\phi(x_{1},x_{2})\in L^{2,2}(\mathbf{R}\times \mathbf{R}^{d}),$ the Fourier transform is defined as
\begin{equation*}
\hat{\phi}(\xi,\tilde{\xi})=\int_{\mathbf{R}}\int_{\mathbf{R}^{d}}\phi(x_{1},x_{2})e^{-i(x_{1},x_{2})\cdot(\xi,\tilde{\xi})}dx_{1}dx_{2}.
\end{equation*}
For $\Phi=(\phi_{1},\cdots,\phi_{r})^{\mathrm{T}}$ and $\Psi=(\psi_{1},\cdots,\psi_{s})^{\mathrm{T}}$, if $\hat{\phi_{i}}(\xi,\tilde{\xi})\overline{\hat{\psi}_{i^{'}}(\xi,\tilde{\xi})}$ is integrable for any $1\leq i \leq r$ and $1\leq i^{'} \leq s $, define a $r \times s$ matrix
\begin{equation*}
[\hat{\Phi},\hat{\Psi}](\xi,\tilde{\xi})=\bigg(\sum\limits_{j_{1}\in \mathbf{Z}}\sum\limits_{j_{2}\in \mathbf{Z}^{d}}\hat{\phi_{i}}(\xi+2j_{1}\pi,\tilde{\xi}+2j_{2}\pi)\overline{\hat{\psi}_{i^{'}}(\xi+2j_{1}\pi,\tilde{\xi}+2j_{2}\pi)}\bigg)_{1\leq i \leq r,\; 1\leq i^{'} \leq s }.
\end{equation*}
\begin{Lm}\label{Lm2.7}
Suppose that $\Phi,\Psi\in \mathcal{L}^{2,2}(\mathbf{R}\times\mathbf{R}^{d})$, then $[\hat{\Phi},\hat{\Psi}](\xi,\tilde{\xi})$ is continuous.
\end{Lm}
{\bfseries Proof} Since $\phi_{i},\psi_{i^{'}}\in \mathcal{L}^{2,2}\subset\mathcal{L}^{1,1}\subset L^{1,1}(\mathbf{R}\times\mathbf{R}^{d}),$ then $\hat{\phi}_{i}$ and $\hat{\psi}_{i^{'}}$ are continuous.

Let $F(\xi,\tilde{\xi})=\hat{\phi}_{i}(\xi,\tilde{\xi})\overline{\hat{\psi}_{i^{'}}(\xi,\tilde{\xi})}$. Then it follows from Poisson's summation formula  that
\begin{eqnarray*}
&&\bigg|\sum\limits_{j_{1}\in \mathbf{Z}}\sum\limits_{j_{2}\in \mathbf{Z}^{d}}\hat{\phi}_{i}(\xi+2j_{1}\pi,\tilde{\xi}+2j_{2}\pi)\overline{\hat{\psi}_{i^{'}}(\xi+2j_{1}\pi,\tilde{\xi}+2j_{2}\pi)}\bigg|\\
&&=\bigg|\sum\limits_{j_{1}\in \mathbf{Z}}\sum\limits_{j_{2}\in \mathbf{Z}^{d}}F(\xi+2j_{1}\pi,\tilde{\xi}+2j_{2}\pi)\bigg|\\
&&=\bigg|\frac{1}{2\pi}\sum\limits_{j_{1}\in \mathbf{Z}}\sum\limits_{j_{2}\in \mathbf{Z}^{d}}\check{F}(j_{1},j_{2})e^{-i(j_{1},j_{2})\cdot(\xi,\tilde{\xi})}\bigg|\\
&&\leq\frac{1}{2\pi}\sum\limits_{j_{1}\in \mathbf{Z}}\sum\limits_{j_{2}\in \mathbf{Z}^{d}}|\phi_{i}\ast\overline{\psi_{i^{'}}}(j_{1},j_{2})|\\
&&\leq\frac{1}{2\pi}\sum\limits_{j_{1}\in \mathbf{Z}}\sum\limits_{j_{2}\in \mathbf{Z}^{d}}\int_{\mathbf{R}}\int_{\mathbf{R}^{d}}|\phi_{i}(x_{1},x_{2})|\;|\psi_{i^{'}}(j_{1}-x_{1},j_{2}-x_{2})|dx_{1}dx_{2}\\
&&=\frac{1}{2\pi}\sum\limits_{j_{1}\in \mathbf{Z}}\sum\limits_{j_{2}\in \mathbf{Z}^{d}}\sum\limits_{k_{1}\in \mathbf{Z}}\sum\limits_{k_{2}\in \mathbf{Z}^{d}}\int_0^1\int_{[0,1]^{d}}|\phi_{i}(x_{1}+k_{1},x_{2}+k_{2})|\;|\psi_{i^{'}}(j_{1}-k_{1}-x_{1},j_{2}-k_{2}-x_{2})|dx_{1}dx_{2}\\
&&=\frac{1}{2\pi}\sum\limits_{k_{1}\in \mathbf{Z}}\sum\limits_{j_{1}\in \mathbf{Z}}\int_0^1\bigg(\int_{[0,1]^{d}}\bigg(\sum\limits_{k_{2}\in \mathbf{Z}^{d}}|\phi_{i}(x_{1}+k_{1},x_{2}+k_{2})|\bigg)\bigg(\sum\limits_{j_{2}\in \mathbf{Z}^{d}}|\psi_{i^{'}}(j_{1}-k_{1}-x_{1},j_{2}-x_{2})|\bigg)dx_{2}\bigg)dx_{1}\\
&&\leq\frac{1}{2\pi}\sum\limits_{k_{1}\in \mathbf{Z}}\sum\limits_{j_{1}\in \mathbf{Z}}\int_0^1\|\phi_{i}(x_{1}+k_{1},\cdot)\|_{\mathcal{L}^{2}}\|\psi_{i^{'}}(j_{1}-k_{1}-x_{1},\cdot)\|_{\mathcal{L}^{2}}dx_{1}\\
&&=\frac{1}{2\pi}\int_0^1\bigg(\sum\limits_{k_{1}\in \mathbf{Z}}\|\phi_{i}(x_{1}+k_{1},\cdot)\|_{\mathcal{L}^{2}}\bigg)\bigg(\sum\limits_{j_{1}\in \mathbf{Z}}\|\psi_{i^{'}}(j_{1}-x_{1},\cdot)\|_{\mathcal{L}^{2}}\bigg)dx_{1}\\
&&\leq\frac{1}{2\pi}\|\phi_{i}\|_{\mathcal{L}^{2,2}}\|\psi_{i^{'}}\|_{\mathcal{L}^{2,2}},
\end{eqnarray*}
which means that
\begin{equation*}
\sum\limits_{j_{1}\in \mathbf{Z}}\sum\limits_{j_{2}\in \mathbf{Z}^{d}}\hat{\phi}_{i}(\xi+2j_{1}\pi,\tilde{\xi}+2j_{2}\pi)\overline{\hat{\psi}_{i^{'}}(\xi+2j_{1}\pi,\tilde{\xi}+2j_{2}\pi)}
\end{equation*}
converges uniformly. Finally, $[\hat{\Phi},\hat{\Psi}](\xi,\tilde{\xi})$ is continuous.
\begin{Lm}\label{L7}
Let $\Phi=(\phi_{1},\cdots,\phi_{r})^{\mathrm{T}}\in \mathcal{L}^{2,2}$. Then the following statements are equivalent to each other.\\
\par $(i)\;\;rank\Big(\hat{\Phi}(\xi+2j_{1}\pi,\tilde{\xi}+2j_{2}\pi)\Big)_{(j_{1},j_{2})\in \mathbf{Z}\times\mathbf{Z}^{d}}$ is a constant function on $\mathbf{R}\times\mathbf{R}^{d}$.\\
\par $(ii)\;\;rank[\hat{\Phi},\hat{\Phi}](\xi,\tilde{\xi})$ is a constant function on $\mathbf{R}\times\mathbf{R}^{d}$.\\
\par $(iii)$ There exists a positive constant C independent of $\xi$ and $\tilde{\xi}$ such that\\
\begin{equation*}
C^{-1}[\hat{\Phi},\hat{\Phi}](\xi,\tilde{\xi})
\leq [\hat{\Phi},\hat{\Phi}](\xi,\tilde{\xi})\overline{[\hat{\Phi},\hat{\Phi}](\xi,\tilde{\xi})^{\mathrm{T}}}\leq C[\hat{\Phi},\hat{\Phi}](\xi,\tilde{\xi}),\ (\xi,\tilde{\xi})\in [-\pi,\pi]\times[-\pi,\pi]^{d}.
\end{equation*}
\end{Lm}
{\bfseries Proof}    Let $G=\Big(\hat{\Phi}(\xi+2j_{1}\pi,\tilde{\xi}+2j_{2}\pi)\Big)_{(j_{1},j_{2})\in \mathbf{Z}\times\mathbf{Z}^{d}}$. Note that $G\overline{G^{\mathrm{T}}}=[\hat{\Phi},\hat{\Phi}](\xi,\tilde{\xi})$. Then the equivalence of (i) and (ii) holds.
Now, we begin to prove the equivalence of (ii) and (iii).

 Note that $[\hat{\Phi},\hat{\Phi}](\xi,\tilde{\xi})$ is a positive semi-definite  Hermite matrix. Then all eigenvalues of $[\hat{\Phi},\hat{\Phi}](\xi,\tilde{\xi})$ satisfy $\lambda_{i}(\xi,\tilde{\xi})\geq 0,\;i=1,\cdots,r$, which are assumed to be ordered as $\lambda_{1}(\xi,\tilde{\xi})\geq \lambda_{2}(\xi,\tilde{\xi})\geq\cdots\geq \lambda_{r}(\xi,\tilde{\xi}).$ Moreover, there exist $r\times r$ matrices $A(\xi,\tilde{\xi})$ such that $\overline{A(\xi,\tilde{\xi})^{\mathrm{T}}}A(\xi,\tilde{\xi})=I_{r}$ and
\begin{equation}\label{S1}
\overline{A(\xi,\tilde{\xi})^{\mathrm{T}}}[\hat{\Phi},\hat{\Phi}](\xi,\tilde{\xi})A(\xi,\tilde{\xi})=diag\bigg(\lambda_{1}(\xi,\tilde{\xi}),\cdots,\lambda_{r}(\xi,\tilde{\xi})\bigg).
\end{equation}
Since $[\hat{\Phi},\hat{\Phi}](\xi,\tilde{\xi})$ is continuous and $2\pi$-periodic, then $\lambda_{i}(\xi,\tilde{\xi})$ are continuous and $2\pi$-periodic for all $1\leq i\leq r$.
\par Let $k_{1}(\xi,\tilde{\xi})=rank[\hat{\Phi},\hat{\Phi}](\xi,\tilde{\xi}).$ If $(ii)$ holds, that is, $k_{1}(\xi,\tilde{\xi})$ is a constant $k_{1}$. Then $\lambda_{i}(\xi,\tilde{\xi})>0$ for all $(\xi,\tilde{\xi})\in \mathbf{R}\times\mathbf{R}^{d}$ and $1\leq i \leq k_{1}$ and\\
\begin{equation}\label{S2}
\lambda_{i}(\xi,\tilde{\xi})\equiv 0\;\;\;for\;(\xi,\tilde{\xi})\in \mathbf{R}\times\mathbf{R}^{d}\;\;and \;\;k_{1}+1\leq i\leq r.
\end{equation}
Moreover, it follows from the continuity and periodicity of $\lambda_{i}(\xi,\tilde{\xi})$ that there exists a positive constant C such that\\
\begin{equation}\label{S3}
C^{-1}\leq \lambda_{i}(\xi,\tilde{\xi})\leq C\;\;\;for\;(\xi,\tilde{\xi})\in \mathbf{R}\times\mathbf{R}^{d}\;\;and\;\; 1\leq i\leq k_{1}.
\end{equation}
This together with $(\ref{S1})$ and $(\ref{S2})$ obtains $(iii)$.
\par If $(iii)$ holds, then it follows from $(\ref{S1})$ that
\begin{equation*}
C^{-1}\lambda_{i}(\xi,\tilde{\xi})\leq \lambda_{i}(\xi,\tilde{\xi})^{2}\leq C\lambda_{i}(\xi,\tilde{\xi})\;\;\;\;for\;\; 1\leq i\leq r\;\;\;and \;\;(\xi,\tilde{\xi})\in \mathbf{R}\times\mathbf{R}^{d}.
\end{equation*}
Thus either $\lambda_{i}(\xi,\tilde{\xi})=0$ or $C^{-1}\leq \lambda_{i}(\xi,\tilde{\xi} )\leq  C$. Hence, rank$[\hat{\Phi},\hat{\Phi}](\xi,\tilde{\xi})$ is a constant by the continuity. This completes the proof of $(ii)\Leftrightarrow (iii)$.
\begin{Lm}\label{L8}
Let $\Phi\in \mathcal{L}^{2,2}$ satisfy $rank\Big(\hat{\Phi}(\xi+2k_{1}\pi,\tilde{\xi}+2k_{2}\pi)\Big)_{(k_{1},k_{2})\in \mathbf{Z}\times \mathbf{Z}^{d}}=k_{0}\geq 1$ for all $(\xi,\tilde{\xi})\in \mathbf{R}\times\mathbf{R}^{d}$. Then there exists a finite  set $\Lambda=\{(\eta_{\lambda_{1}},\eta_{\lambda_{2}})\}\subset[-\pi,\pi]\times[-\pi,\pi]^{d}$, $r\times r$ $2\pi$-periodic nonsingular matrices $P_{(\eta_{\lambda_{1}},\eta_{\lambda_{2}})}(\xi,\tilde{\xi})$ and $(K_{\eta_{\lambda_{1}}},K_{\eta_{\lambda_{2}}}) \subset\mathbf{Z}\times \mathbf{Z}^{d}$ with cardinality $\#(K_{\eta_{\lambda_{1}}},K_{\eta_{\lambda_{2}}})=k_{0},\;2\pi$-periodic $C^{\infty}$ functions $h_{(\eta_{\lambda_{1}},\eta_{\lambda_{2}})}(\xi,\tilde{\xi})$ on $\mathbf{R}\times\mathbf{R}^{d}$, which satisfy the following properties:
\par (i) There exist $\Psi_{1,(\eta_{\lambda_{1}},\eta_{\lambda_{2}})}$ and $\Psi_{2,(\eta_{\lambda_{1}},\eta_{\lambda_{2}})}$ which are functions from $\mathbf{R}^{d+1}$ to $C^{k_{0}}$ and $C^{r-k_{0}}$ respectively, defined by
\begin{equation}
P_{(\eta_{\lambda_{1}},\eta_{\lambda_{2}})}(\xi,\tilde{\xi})\hat{\Phi}(\xi,\tilde{\xi})
=\left(
\begin{array}{c} 
\hat{\Psi}_{1,(\eta_{\lambda_{1}},\eta_{\lambda_{2}})}(\xi,\tilde{\xi})\\ 
\hat{\Psi}_{2,(\eta_{\lambda_{1}},\eta_{\lambda_{2}})}(\xi,\tilde{\xi})\\ 
\end{array}
\right)\;\;\;for\;(\xi,\tilde{\xi})\in\mathbf{R}\times\mathbf{R}^{d}.
\end{equation}
Moreover, for sufficiently small $\delta_{(\eta_{\lambda_{1}},\eta_{\lambda_{2}})}>0$, one has
\begin{equation}\label{S5}
rank\bigg(\hat{\Psi}_{1,(\eta_{\lambda_{1}},\eta_{\lambda_{2}})}(\xi+2k_{1}\pi,\tilde{\xi}+2k_{2}\pi)\bigg)_{(k_{1},k_{2})\in (K_{\eta_{\lambda_{1}}}, K_{\eta_{\lambda_{2}}})}=k_{0},\ (\xi,\tilde{\xi})\in B\Big((\eta_{\lambda_{1}},\eta_{\lambda_{2}}),2\delta_{(\eta_{\lambda_{1}},\eta_{\lambda_{2}})}\Big)
\end{equation}
and
\begin{equation}\label{S6}
\hat{\Psi}_{2,(\eta_{\lambda_{1}},\eta_{\lambda_{2}})}(\xi,\tilde{\xi})=0,\ (\xi,\tilde{\xi})\in B\Big(({\eta_{\lambda_{1}},\eta_{\lambda_{2}}}),2\delta_{(\eta_{\lambda_{1}},\eta_{\lambda_{2}})}\Big)+2\pi (\mathbf{Z}\times\mathbf{Z}^{d}).
\end{equation}
\par (ii)
\begin{equation}\label{r13}
\sum\limits_{(\eta_{\lambda_{1}},\eta_{\lambda_{2}})\in \Lambda}h_{(\eta_{\lambda_{1}},\eta_{\lambda_{2}})}(\xi,\tilde{\xi})=1,\ (\xi,\tilde{\xi})\in \mathbf{R}\times\mathbf{R}^{d}
\end{equation}
and
\begin{equation}\label{r14}
supp\;h_{(\eta_{\lambda_{1}},\eta_{\lambda_{2}})}(\xi,\tilde{\xi})\subset B\Big((\eta_{\lambda_{1}},\eta_{\lambda_{2}}),\delta_{(\eta_{\lambda_{1}},\eta_{\lambda_{2}})}\Big)+2\pi (\mathbf{Z}\times\mathbf{Z}^{d}).
\end{equation}
\end{Lm}
{\bfseries Proof}  For any $(\eta_{1},\eta_{2})\in[-\pi,\pi]\times [-\pi,\pi]^{d}$, there exists a $r\times r$ nonsingular matrix $P_{(\eta_{1},\eta_{2})}$, a $k_{0}\times k_{0}$ nonsingular matrix $A_{(\eta_{1},\eta_{2})}$ and $(K_{\eta_{1}},K_{\eta_{2}})\subset \mathbf{Z}\times\mathbf{Z}^{d}$ with cardinality $k_{0}$ such that
\begin{equation*}
P_{(\eta_{1},\eta_{2})}\bigg(\hat{\Phi}(\eta_{1}+2k_{1}\pi,\eta_{2}+2k_{2}\pi)\bigg)_{(k_{1},k_{2})\in( K_{\eta_{1}},K_{\eta_{2}})}
=\left(
\begin{array}{c} 
A_{(\eta_{1},\eta_{2})}\\ 
0\\ 
\end{array}
\right).
\end{equation*}
Write
\begin{equation*}
P_{(\eta_{1},\eta_{2})}\bigg(\hat{\Phi}(\xi+2k_{1}\pi,\tilde{\xi}+2k_{2}\pi)\bigg)_{(k_{1},k_{2})\in (K_{\eta_{1}}, K_{\eta_{2}})}
=\left(
\begin{array}{c} 
A_{(\eta_{1},\eta_{2})}+R_{1}(\xi,\tilde{\xi})\\ 
R_{2}(\xi,\tilde{\xi})\\ 
\end{array}
\right).
\end{equation*}
By the continuity of $\hat{\Phi},\;R_{1}(\xi,\tilde{\xi})$ and $R_{2}(\xi,\tilde{\xi})$ are continuous, $R_{1}(\eta_{1},\eta_{2})= R_{2}(\eta_{1},\eta_{2})=0$. Thus $$\sup\limits_{(\xi,\tilde{\xi})\in B((\eta_{1},\eta_{2}),5\delta_{(\eta_{1},\eta_{2})})}\|R_{1}(\xi,\tilde{\xi})\|+\|R_{2}(\xi,\tilde{\xi})\|$$ is sufficiently small for any sufficiently small  $\delta_{(\eta_{1},\eta_{2})}>0$. Let $H(x,y)$ be a nonnegative $C^{\infty}$ function on $\mathbf{R}^{d+1}$ such that\\
\begin{equation}\label{S7}
H(x,y)=
\begin{cases}
1,   \;\;\;\;\;&\mbox{if $|(x,y)|\leq 1$},\\
0,   \;\;\;\;\;&\mbox{if $|(x,y)|\geq 5/4$ }.
\end{cases}
\end{equation}
Then $A_{(\eta_{1},\eta_{2})}+H\Big((\xi-\eta_{1})/4\delta_{(\eta_{1},\eta_{2})},(\tilde{\xi}-\eta_{2})/4\delta_{(\eta_{1},\eta_{2})}\Big)R_{1}(\xi,\tilde{\xi})$ is a $k_{0}\times k_{0}$ nonsingular matrix for all $(\xi,\tilde{\xi})\in \mathbf{R}\times\mathbf{R}^{d}$ when $\delta_{(\eta_{1},\eta_{2})}$ is chosen sufficiently small. For $(\xi,\tilde{\xi})\in\mathbf{R}\times\mathbf{R}^{d}$, set
\begin{equation*}
\alpha_{(\eta_{1},\eta_{2})}(\xi,\tilde{\xi})=A_{(\eta_{1},\eta_{2})}+\sum\limits_{j_{1}\in \mathbf{Z}}\sum\limits_{j_{2}\in \mathbf{Z}^{d}}H\bigg(\frac{\xi+2j_{1}\pi-\eta_{1}}{4\delta_{(\eta_{1},\eta_{2})}},\frac{\tilde{\xi}+2j_{2}\pi-\eta_{2}}{4\delta_{(\eta_{1},\eta_{2})}}\bigg)R_{1}(\xi+2j_{1}\pi,\tilde{\xi}+2j_{2}\pi)
\end{equation*}
and
\begin{equation*}
\beta_{(\eta_{1},\eta_{2})}(\xi,\tilde{\xi})=\sum\limits_{j_{1}\in \mathbf{Z}}\sum\limits_{j_{2}\in \mathbf{Z}^{d}}H\bigg(\frac{\xi+2j_{1}\pi-\eta_{1}}{2\delta_{(\eta_{1},\eta_{2})}},\frac{\tilde{\xi}+2j_{2}\pi-\eta_{2}}{2\delta_{(\eta_{1},\eta_{2})}}\bigg)R_{2}(\xi+2j_{1}\pi,\tilde{\xi}+2j_{2}\pi).
\end{equation*}
Then $\alpha_{(\eta_{1},\eta_{2})}(\xi,\tilde{\xi})$ is $2\pi$-periodic and nonsingular when $\delta_{(\eta_{1},\eta_{2})}$ is chosen sufficiently small.
Let
\begin{equation*}
P_{(\eta_{1},\eta_{2})}(\xi,\tilde{\xi})=P_{(\eta_{1},\eta_{2})}+
\begin{pmatrix}
0&0\\
-\beta_{(\eta_{1},\eta_{2})}(\xi,\tilde{\xi})\Big(\alpha_{(\eta_{1},\eta_{2})}(\xi,\tilde{\xi})\Big)^{-1}&0
\end{pmatrix}
P_{(\eta_{1},\eta_{2})},\;\;\;(\xi,\tilde{\xi})\in \mathbf{R}\times\mathbf{R}^{d}.
\end{equation*}
Then $P_{(\eta_{1},\eta_{2})}(\xi,\tilde{\xi})$ is a $r\times r$ $2\pi$-periodic nonsingular matrix for any $(\xi,\tilde{\xi})\in \mathbf{R}\times\mathbf{R}^{d}$. Note that $\alpha_{(\eta_{1},\eta_{2})}(\xi,\tilde{\xi})=A_{(\eta_{1},\eta_{2})}+R_{1}(\xi,\tilde{\xi})$ and $\beta_{(\eta_{1},\eta_{2})}(\xi,\tilde{\xi})=R_{2}(\xi,\tilde{\xi})$ for $(\xi,\tilde{\xi})\in B\Big((\eta_{1},\eta_{2}),2\delta_{(\eta_{1},\eta_{2})}\Big)$. Thus for any $(\xi,\tilde{\xi})\in B\Big((\eta_{1},\eta_{2}),2\delta_{(\eta_{1},\eta_{2})}\Big)$, one has
\begin{equation}\label{S8}
P_{(\eta_{1},\eta_{2})}(\xi,\tilde{\xi})\bigg(\hat{\Phi}(\xi+2k_{1}\pi,\tilde{\xi}+2k_{2}\pi)\bigg)_{(k_{1},k_{2})\in (K_{\eta_{1}},K_{\eta_{2}})}=
\begin{pmatrix}
A_{(\eta_{1},\eta_{2})}+R_{1}(\xi,\tilde{\xi})\\
0
\end{pmatrix}.
\end{equation}
Define $\Psi_{1,(\eta_{1},\eta_{2})}=(\psi_{1,(\eta_{1},\eta_{2}),1},\ldots,\psi_{1,(\eta_{1},\eta_{2}),k_{0}})^{\mathrm{T}}$ and $\Psi_{2,(\eta_{1},\eta_{2})}=(\psi_{2,(\eta_{1},\eta_{2}),1},\ldots,\psi_{2,(\eta_{1},\eta_{2}),r-k_{0}})^{\mathrm{T}}$ by\\
\begin{equation}\label{S9}
\begin{pmatrix}
\hat{\Psi}_{1,(\eta_{1},\eta_{2})}(\xi,\tilde{\xi})\\
\hat{\Psi}_{2,(\eta_{1},\eta_{2})}(\xi,\tilde{\xi})
\end{pmatrix}
=P_{(\eta_{1},\eta_{2})}(\xi,\tilde{\xi})\hat{\Phi}(\xi,\tilde{\xi}),\;\;\;\;(\xi,\tilde{\xi})\in \mathbf{R}\times\mathbf{R}^{d}.
\end{equation}
Recall that the Fourier coefficients of all entries of $P_{(\eta_{1},\eta_{2})}(\xi,\tilde{\xi})$  belong to  $\ell^{1,1}$. Thus by Lemma $\ref{lm 2.4}$, we obtain $\Psi_{1,(\eta_{1},\eta_{2})}\in\mathcal{L}^{1,1}$ and $\Psi_{2,(\eta_{1},\eta_{2})}\in\mathcal{L}^{2,2}$. By $(\ref{S8})$ and $(\ref{S9})$, we obtain
\begin{equation}\label{S10}
\hat{\Psi}_{2,(\eta_{1},\eta_{2})}(\xi+2k_{1}\pi,\tilde{\xi}+2k_{2}\pi)=0\;\;\;for \;(\xi,\tilde{\xi})\in B\Big((\eta_{1},\eta_{2}),2\delta_{(\eta_{1},\eta_{2})}\Big)\;and\;(k_{1},k_{2})\in (K_{\eta_{1}},K_{\eta_{2}}).
\end{equation}
For $(\xi,\tilde{\xi})\in B\Big((\eta_{1},\eta_{2}),2\delta_{(\eta_{1},\eta_{2})}\Big)$, we have $$rank\bigg(\hat{\Psi}_{1,(\eta_{1},\eta_{2})}(\xi+2k_{1}\pi,\tilde{\xi}+2k_{2}\pi)\bigg)_{(k_{1},k_{2})\in(K_{\eta_{1}},K_{\eta_{2}})}=rank\Big(A_{(\eta_{1},\eta_{2})}+R_{1}(\xi,\tilde{\xi})\Big)=k_{0}$$ and
\begin{eqnarray*}
&&rank
\begin{pmatrix}
A_{(\eta_{1},\eta_{2})}+R_{1}(\xi,\tilde{\xi})&\hat{\Psi}_{1,(\eta_{1},\eta_{2})}(\xi+2k_{1}^{'}\pi,\tilde{\xi}+2k_{2}^{'}\pi)\\
0&\hat{\Psi}_{2,(\eta_{1},\eta_{2})}(\xi+2k_{1}^{'}\pi,\tilde{\xi}+2k_{2}^{'}\pi)
\end{pmatrix}
_{(k_{1}^{'},k_{2}^{'})\in(\mathbf{Z}\times\mathbf{Z}^{d})\backslash (K_{\eta_{1}},K_{\eta_{2}})}\\
&&= rank
\begin{pmatrix}
\hat{\Psi}_{1,(\eta_{1},\eta_{2})}(\xi+2k_{1}\pi,\tilde{\xi}+2k_{2}\pi)\\
\hat{\Psi}_{2,(\eta_{1},\eta_{2})}(\xi+2k_{1}\pi,\tilde{\xi}+2k_{2}\pi)
\end{pmatrix}
_{(k_{1},k_{2})\in\mathbf{Z}\times\mathbf{Z}^{d}}\\
&&=rank\bigg(\hat{\Phi}(\xi+2k_{1}\pi,\tilde{\xi}+2k_{2}\pi)\bigg)_{(k_{1},k_{2})\in\mathbf{Z}\times\mathbf{Z}^{d}}=k_{0}.
\end{eqnarray*}
Thus $\hat{\Psi}_{2,(\eta_{1},\eta_{2})}(\xi+2k_{1}^{'}\pi,\tilde{\xi}+2k_{2}^{'}\pi)=0$ for $(\xi,\tilde{\xi})\in B\Big((\eta_{1},\eta_{2}),2\delta_{(\eta_{1},\eta_{2})}\Big)$ and $(k_{1}^{'},k_{2}^{'})\in(\mathbf{Z}\times\mathbf{Z}^{d})\backslash (K_{\eta_{1}},K_{\eta_{2}})$. This together with $(\ref{S10})$ obtains\\
\begin{equation*}
\hat{\Psi}_{2,(\eta_{1},\eta_{2})}(\xi,\tilde{\xi})=0\;\;\;for\;(\xi,\tilde{\xi})\in B\Big((\eta_{1},\eta_{2}),2\delta_{(\eta_{1},\eta_{2})}\Big)+2\pi(\mathbf{Z}\times\mathbf{Z}^{d}).
\end{equation*}
\par For the family $\Big\{B\Big((\eta_{1},\eta_{2}),\delta_{(\eta_{1},\eta_{2})}/2\Big):(\eta_{1},\eta_{2})\in[-\pi,\pi]\times[-\pi,\pi]^{d}\Big\}$, it follows from finite covering theorem that there exists a finite set $\Lambda\subset[-\pi,\pi]\times[-\pi,\pi]^{d}$ such that
\begin{equation*}
\bigcup\limits_{(\eta_{\lambda_{1}},\eta_{\lambda_{2}})\in \Lambda}B\Big((\eta_{\lambda_{1}},\eta_{\lambda_{2}}),\delta_{(\eta_{\lambda_{1}},\eta_{\lambda_{2}})}/2\Big)\supset [-\pi,\pi]\times [-\pi,\pi]^{d}.
\end{equation*}
Then, there exist $2\pi$-periodic $C^{\infty}$ functions $h_{(\eta_{\lambda_{1}},\eta_{\lambda_{2}})}(\xi,\tilde{\xi})$ such that $(\ref{r13})$ and $(\ref{r14})$ hold.
\begin{Lm}\label{L9}
Let $\phi\in\mathcal{L}^{p,q}$ if $1\leq p,q<\infty$ and $\phi\in \mathcal{W}(L^{1,1})$ if $p=\infty$ or $q=\infty$. Assume that $\sum\limits_{j_{1}\in\mathbf{Z}}\sum\limits_{j_{2}\in\mathbf{Z}^{d}}\phi(\cdot-j_{1},\cdot-j_{2})=0$. Then for $0< \varepsilon_{1},\;\varepsilon_{2}< 1$ and any function $h$ on $\mathbf{R}\times\mathbf{R}^{d}$ satisfying
\begin{equation}\label{S11}
|h(x_{1},y_{1})-h(x_{2},y_{2})|\leq C|x_{1}-x_{2}|\;|y_{1}-y_{2}|\Big(1+min(|x_{1}|,|x_{2}|)\Big)^{-1-\varepsilon_{1}}\Big(1+min(|y_{1}|,|y_{2}|)\Big)^{-d-\varepsilon_{2}},
\end{equation}
we obtain
\begin{equation}\label{rr2.15}
\lim\limits_{n\rightarrow\infty}2^{-n(d+1)}\Big\|\sum\limits_{j_{1}\in\mathbf{Z}}\sum\limits_{j_{2}\in\mathbf{Z}^{d}}h(2^{-n}j_{1},2^{-n}j_{2})\phi(\cdot-j_{1},\cdot-j_{2})\Big\|_{\mathcal{L}^{p,q}}=0.
\end{equation}
\end{Lm}
{\bfseries Proof} For $m\in\mathbf{Z}^{d}$, let $|m|=max\{|m_{1}|,|m_{2}|,\cdots,|m_{d}|\}$. If $1\leq p,q<\infty$, then for any $\varepsilon > 0$, there exists a positive integer $N$ such that\\
\begin{equation}\label{r20}
\bigg\|\sum\limits_{|j_{1}|\geq N}\bigg(\int_{[0,1]^{d}}\bigg(\sum\limits_{j_{2}\in\mathbf{Z}^{d}} \Big|\phi(x_{1}+j_{1},x_{2}+j_{2})\Big|\bigg)^{q}dx_{2}\bigg)^{1/q}\bigg\|_{L^{p}([0,1])}\leq \frac{\varepsilon}{4}
\end{equation}
and
\begin{equation}\label{rr20}
\bigg\|\sum\limits_{|j_{1}|< N}\bigg(\int_{[0,1]^{d}}\bigg(\sum\limits_{|j_{2}|\geq N}\Big|\phi(x_{1}+j_{1},x_{2}+j_{2})\Big|\bigg)^{q}dx_{2}\bigg)^{1/q}\bigg\|_{L^{p}([0,1])}\leq \frac{\varepsilon}{4}.
\end{equation}
Set
\begin{eqnarray*}
\phi_{1}(x_{1},x_{2})=&&\phi(x_{1},x_{2})\chi_{O_{N_{0}}}(x_{1},x_{2})\\
&&+\bigg(\sum\limits_{|j_{1}|\geq N}\sum\limits_{|j_{2}|\geq N}+\sum\limits_{|j_{1}|\geq N}\sum\limits_{|j_{2}|< N}+\sum\limits_{|j_{1}|< N}\sum\limits_{|j_{2}|\geq N}\bigg)\phi(x_{1}+j_{1},x_{2}+j_{2})\chi_{[0,1]\times[0,1]^{d}}(x_{1},x_{2}),
\end{eqnarray*}
where $O_{N_{0}}=\bigcup\limits_{\{|j_{1}|< N,\;|j_{2}|< N\}}\Big((j_{1},j_{2})+[0,1]\times[0,1]^{d}\Big)$.
Then for $(x_{1},x_{2})\in[0,1]\times[0,1]^{d}$,
\begin{equation}\label{r21}
\sum_{j_{1}\in\mathbf{Z}}\sum_{j_{2}\in\mathbf{Z}^{d}}\phi_{1}(x_{1}-j_{1},x_{2}-j_{2})=:I_{1}+I_{2}.
\end{equation}
It is easy to verify that
\begin{eqnarray}\label{L22}
I_{1}=&&\sum\limits_{j_{1}\in\mathbf{Z}}\sum_{j_{2}\in\mathbf{Z}^{d}}\phi(x_{1}-j_{1},x_{2}-j_{2})\chi_{O_{N_{0}}}(x_{1}-j_{1},x_{2}-j_{2})\nonumber\\
=&&\sum\limits_{|j_{1}|< N}\sum\limits_{|j_{2}|< N}\phi(x_{1}-j_{1},x_{2}-j_{2}),\\
I_{2}=&&\sum_{j_{1}\in\mathbf{Z}}\sum\limits_{j_{2}\in\mathbf{Z}^{d}}\Big(\sum\limits_{|k_{1}|\geq N}\sum\limits_{|k_{2}|\geq N}+\sum\limits_{|k_{1}|\geq N}\sum\limits_{|k_{2}|< N}+\sum\limits_{|k_{1}|< N}\sum\limits_{|k_{2}|\geq N}\Big)\nonumber\\
&&\times\phi(x_{1}+k_{1}-j_{1},x_{2}+k_{2}-j_{2})\chi_{[0,1]\times[0,1]^{d}}(x_{1}-j_{1},x_{2}-j_{2})\nonumber\\
=&&\Big(\sum\limits_{|k_{1}|\geq N}\sum\limits_{|k_{2}|\geq N}+\sum\limits_{|k_{1}|\geq N}\sum\limits_{|k_{2}|< N}+\sum\limits_{|k_{1}|< N}\sum\limits_{|k_{2}|\geq N}\Big)\phi(x_{1}+k_{1},x_{2}+k_{2}).\nonumber
\end{eqnarray}
This together with $(\ref{r21})$ and $(\ref{L22})$ obtains
\begin{equation}\label{r2.23}
\sum_{j_{1}\in\mathbf{Z}}\sum_{j_{2}\in\mathbf{Z}^{d}}\phi_{1}(x_{1}-j_{1},x_{2}-j_{2})=\sum_{j_{1}\in\mathbf{Z}}\sum_{j_{2}\in\mathbf{Z}^{d}}\phi(x_{1}-j_{1},x_{2}-j_{2})=0.
\end{equation}
Moreover, we obtain
\begin{eqnarray*}
\|\phi_{1}-\phi\|_{\mathcal{L}^{p,q}}=&&\bigg\|\sum\limits_{j_{1}\in\mathbf{Z}}\bigg(\int_{[0,1]^{d}}\bigg(\sum\limits_{j_{2}\in\mathbf{Z}^{d}}\bigg|\phi(x_{1}-j_{1},x_{2}-j_{2})\chi_{O_{N_{0}}}(x_{1}-j_{1},x_{2}-j_{2})\\
&&+\bigg(\sum\limits_{|k_{1}|\geq N}\sum\limits_{|k_{2}|\geq N}+\sum\limits_{|k_{1}|\geq N}\sum\limits_{|k_{2}|< N}+\sum\limits_{|k_{1}|< N}\sum\limits_{|k_{2}|\geq N}\bigg)\phi(x_{1}+k_{1}-j_{1},x_{2}+k_{2}-j_{2})\\
&&\times\chi_{[0,1]\times[0,1]^{d}}(x_{1}-j_{1},x_{2}-j_{2})-\phi(x_{1}-j_{1},x_{2}-j_{2})\bigg|\bigg)^{q}dx_{2}\bigg)^{1/q}\bigg\|_{L^{p}([0,1])}\\
\leq&& E_{1}+E_{2}.
\end{eqnarray*}
Next, we discuss $E_{1}$ and $E_{2}$ respectively. In fact, one has
\begin{eqnarray}\label{r23}
E_{1}=&&\bigg\|\sum\limits_{|j_{1}|\geq N} \bigg(\int_{[0,1]^{d}}\bigg(\sum\limits_{j_{2}\in\mathbf{Z}^{d}}\bigg|\phi(x_{1}-j_{1},x_{2}-j_{2})\chi_{O_{N_{0}}}(x_{1}-j_{1},x_{2}-j_{2})\nonumber\\
&&+\bigg(\sum\limits_{|k_{1}|\geq N}\sum\limits_{|k_{2}|\geq N}+\sum\limits_{|k_{1}|\geq N}\sum\limits_{|k_{2}|< N}+\sum\limits_{|k_{1}|< N}\sum\limits_{|k_{2}|\geq N}\bigg)\phi(x_{1}+k_{1}-j_{1},x_{2}+k_{2}-j_{2})\nonumber\\
&&\times\chi_{[0,1]\times[0,1]^{d}}(x_{1}-j_{1},x_{2}-j_{2})-\phi(x_{1}-j_{1},x_{2}-j_{2})\bigg|\bigg)^{q}dx_{2}\bigg)^{1/q}\bigg\|_{L^{p}([0,1])}\nonumber\\
=&&\bigg\|\sum\limits_{|j_{1}|\geq N} \bigg(\int_{[0,1]^{d}}\bigg(\sum\limits_{j_{2}\in\mathbf{Z}^{d}}\Big|\phi(x_{1}-j_{1},x_{2}-j_{2})\Big|\bigg)^{q}dx_{2}\bigg)^{1/q}\bigg\|_{L^{p}([0,1])}\leq \frac{\varepsilon}{4},
\end{eqnarray}
\begin{eqnarray*}
E_{2}=&&\bigg\|\sum\limits_{|j_{1}|< N} \bigg(\int_{[0,1]^{d}}\bigg(\sum\limits_{j_{2}\in\mathbf{Z}^{d}}\bigg|\phi(x_{1}-j_{1},x_{2}-j_{2})\chi_{O_{N_{0}}}(x_{1}-j_{1},x_{2}-j_{2})\\
&&+\bigg(\sum\limits_{|k_{1}|\geq N}\sum\limits_{|k_{2}|\geq N}+\sum\limits_{|k_{1}|\geq N}\sum\limits_{|k_{2}|< N}+\sum\limits_{|k_{1}|< N}\sum\limits_{|k_{2}|\geq N}\bigg)\phi(x_{1}+k_{1}-j_{1},x_{2}+k_{2}-j_{2})\\
&&\times\chi_{[0,1]\times[0,1]^{d}}(x_{1}-j_{1},x_{2}-j_{2})-\phi(x_{1}-j_{1},x_{2}-j_{2})\bigg|\bigg)^{q}dx_{2}\bigg)^{1/q}\bigg\|_{L^{p}([0,1])}\\
\leq && E_{21}+E_{22},
\end{eqnarray*}
where
\begin{eqnarray*}
E_{21}=&&\bigg\|\sum\limits_{|j_{1}|< N} \bigg(\int_{[0,1]^{d}}\bigg(\sum\limits_{|j_{2}|<N}\bigg|\phi(x_{1}-j_{1},x_{2}-j_{2})\chi_{O_{N_{0}}}(x_{1}-j_{1},x_{2}-j_{2})\nonumber\\
&&+\bigg(\sum\limits_{|k_{1}|\geq N}\sum\limits_{|k_{2}|\geq N}+\sum\limits_{|k_{1}|\geq N}\sum\limits_{|k_{2}|< N}+\sum\limits_{|k_{1}|< N}\sum\limits_{|k_{2}|\geq N}\bigg)\phi(x_{1}+k_{1}-j_{1},x_{2}+k_{2}-j_{2})\nonumber\\
&&\times\chi_{[0,1]\times[0,1]^{d}}(x_{1}-j_{1},x_{2}-j_{2})-\phi(x_{1}-j_{1},x_{2}-j_{2})\bigg|\bigg)^{q}dx_{2}\bigg)^{1/q}\bigg\|_{L^{p}([0,1])}\nonumber\\
\end{eqnarray*}
\begin{eqnarray}\label{r24}
=&&\bigg\|\bigg(\int_{[0,1]^{d}}\bigg|
\bigg(\sum\limits_{|k_{1}|\geq N}\sum\limits_{|k_{2}|\geq N}+\sum\limits_{|k_{1}|\geq N}\sum\limits_{|k_{2}|< N}+\sum\limits_{|k_{1}|< N}\sum\limits_{|k_{2}|\geq N}\bigg)\phi(x_{1}+k_{1},x_{2}+k_{2})\bigg|^{q}dx_{2}\bigg)^{1/q}\bigg\|_{L^{p}([0,1])}\nonumber\\
\leq&&\bigg\|\sum\limits_{|k_{1}|\geq N}\bigg(\int_{[0,1]^{d}}
\bigg(\sum\limits_{k_{2}\in\mathbf{Z}^{d}}\Big|\phi(x_{1}+k_{1},x_{2}+k_{2})\Big|\bigg)^{q}dx_{2}\bigg)^{1/q}\bigg\|_{L^{p}([0,1])}\nonumber\\
&&+\bigg\|\sum\limits_{|k_{1}|< N}\bigg(\int_{[0,1]^{d}}\bigg(\sum\limits_{|k_{2}|\geq N}\Big|\phi(x_{1}+k_{1},x_{2}+k_{2})\Big|\bigg)^{q}dx_{2}\bigg)^{1/q}\bigg\|_{L^{p}([0,1])}\nonumber\\
\leq&&\frac{\varepsilon}{4}+\frac{\varepsilon}{4}=\frac{\varepsilon}{2}
\end{eqnarray}
and
\begin{eqnarray*}
E_{22}=&&\bigg\|\sum\limits_{|j_{1}|< N} \bigg(\int_{[0,1]^{d}}\bigg(\sum\limits_{|j_{2}|\geq N}\bigg|\phi(x_{1}-j_{1},x_{2}-j_{2})\chi_{O_{N_{0}}}(x_{1}-j_{1},x_{2}-j_{2})\\
&&+\bigg(\sum\limits_{|k_{1}|\geq N}\sum\limits_{|k_{2}|\geq N}+\sum\limits_{|k_{1}|\geq N}\sum\limits_{|k_{2}|< N}+\sum\limits_{|k_{1}|< N}\sum\limits_{|k_{2}|\geq N}\bigg)\phi(x_{1}+k_{1}-j_{1},x_{2}+k_{2}-j_{2})\\
&&\times\chi_{[0,1]\times[0,1]^{d}}(x_{1}-j_{1},x_{2}-j_{2})-\phi(x_{1}-j_{1},x_{2}-j_{2})\bigg|\bigg)^{q}dx_{2}\bigg)^{1/q}\bigg\|_{L^{p}([0,1])}\\
=&&\bigg\|\sum\limits_{|j_{1}|< N} \bigg(\int_{[0,1]^{d}}\bigg(\sum\limits_{|j_{2}|\geq N}\Big|\phi(x_{1}-j_{1},x_{2}-j_{2})\Big|\bigg)^{q}dx_{2}\bigg)^{1/q}\bigg\|_{L^{p}([0,1])}\\
\leq&&\frac{\varepsilon}{4}.
\end{eqnarray*}
This together with $(\ref{r23})$ and $(\ref{r24})$ obtains
\begin{equation*}
\|\phi_{1}-\phi\|_{\mathcal{L}^{p,q}}\leq \varepsilon.
\end{equation*}
Therefore, using $(\ref{S11})$, $(\ref{r2.23})$ and the fact that
\begin{equation}\label{r2.25}
supp\;\phi_{1}\subset \bigg\{\Big((j_{1},j_{2})+[0,1]\times[0,1]^{d}\Big):|j_{1}|\leq N,\;|j_{2}|\leq N\bigg\},
\end{equation}
we obtain
\begin{eqnarray*}
&&2^{-n(d+1)p}
\bigg\|\sum\limits_{j_{1}\in\mathbf{Z}}\sum\limits_{j_{2}\in\mathbf{Z}^{d}}h(2^{-n}j_{1},2^{-n}j_{2})\phi_{1}(x_{1}-j_{1}+k_{1},x_{2}-j_{2}+k_{2})\bigg\|_{\mathcal{L}^{p,q}}^{p}\\
=&& 2^{-n(d+1)p}
\int_{[0,1]}\bigg(\sum\limits_{k_{1}\in\mathbf{Z}}\bigg(\int_{[0,1]^{d}}\bigg(\sum\limits_{k_{2}\in\mathbf{Z}^{d}}\Big|\sum\limits_{j_{1}\in\mathbf{Z}}\sum\limits_{j_{2}\in\mathbf{Z}^{d}}\bigg(h(2^{-n}j_{1},2^{-n}j_{2})-h(2^{-n}k_{1},2^{-n}k_{2})\bigg)\\
&&\phi_{1}(x_{1}-j_{1}+k_{1},x_{2}-j_{2}+k_{2})\Big|\bigg)^{q}dx_{2}\bigg)^{1/q}\bigg)^{p}dx_{1}\\
\end{eqnarray*}
\begin{eqnarray*}
\leq&& C^{p}2^{-n(d+1)p}\int_{[0,1]}\bigg(\sum\limits_{k_{1}\in\mathbf{Z}}\bigg(\int_{[0,1]^{d}}\bigg(\sum\limits_{k_{2}\in\mathbf{Z}^{d}}\sum\limits_{j_{1}\in\mathbf{Z}}\sum\limits_{j_{2}\in\mathbf{Z}^{d}}\Big|2^{-n}j_{1}-2^{-n}k_{1}\Big|\;\Big|2^{-n}j_{2}-2^{-n}k_{2}\Big|\\
&&\Big(1+min(|2^{-n}j_{1}|,|2^{-n}k_{1}|)\Big)^{-1-\varepsilon_{1}}\Big(1+min(|2^{-n}j_{2}|,|2^{-n}k_{2}|)\Big)^{-d-\varepsilon_{2}}\\
&&\Big|\phi_{1}(x_{1}-j_{1}+k_{1},x_{2}-j_{2}+k_{2})\Big|\bigg)^{q}dx_{2}\bigg)^{1/q}\bigg)^{p}dx_{1}\\
\leq&& C_{1}^{p}(N)C^{p}2^{-n(d+1)p}2^{-2np}\int_{[0,1]}\bigg(\sum\limits_{k_{1}\in\mathbf{Z}}\bigg(\int_{[0,1]^{d}}\bigg(\sum\limits_{k_{2}\in\mathbf{Z}^{d}}\sum\limits_{|j_{1}-k_{1}|<N}\sum\limits_{|j_{2}-k_{2}|<N}\Big(1+2^{-n}|k_{1}|\Big)^{-1-\varepsilon_{1}}\\
&&\Big(1+2^{-n}|k_{2}|\Big)^{-d-\varepsilon_{2}}\Big|\phi_{1}(x_{1}-j_{1}+k_{1},x_{2}-j_{2}+k_{2})\Big|\bigg)^{q}dx_{2}\bigg)^{1/q}\bigg)^{p}dx_{1}\\
\leq&&
C_{1}^{p}(N)C^{p}2^{-n(d+1)p}2^{-2np}\int_{[0,1]}\bigg(\sum\limits_{k_{1}\in\mathbf{Z}}\Big(1+2^{-n}|k_{1}|\Big)^{-1-\varepsilon_{1}}\bigg(\int_{[0,1]^{d}}\bigg(\sum\limits_{k_{2}\in\mathbf{Z}^{d}}\Big(1+2^{-n}|k_{2}|\Big)^{-d-\varepsilon_{2}}\\
&&\sum\limits_{|j_{1}-k_{1}|<N}\sum\limits_{|j_{2}-k_{2}|<N}\Big|\phi_{1}(x_{1}-j_{1}+k_{1},x_{2}-j_{2}+k_{2})\Big|\bigg)^{q}dx_{2}\bigg)^{1/q}\bigg)^{p}dx_{1}\\
\leq&&2^{n(1+\varepsilon_{1})p}2^{n(d+\varepsilon_{2})p}2^{-n(d+1)p}2^{-2np}C^{p}C^{p}(N)\|\phi_{1}\|_{\mathcal{L}^{p,q}}^{p}\\
\leq&&2^{n(\varepsilon_{1}-1)p}2^{n(\varepsilon_{2}-1)p}C^{p}C^{p}(N)\Big(\|\phi\|_{\mathcal{L}^{p,q}}+\varepsilon\Big)^{p},
\end{eqnarray*}
where $C(N)$ is a positive constant depending only on $N,\;d,\;\varepsilon_{1}$ and $\varepsilon_{2}$,  $C$ is the constant in $(\ref{S11})$. Then, $(\ref{rr2.15})$ holds. 
\par If $p=\infty$ or $q=\infty$, then for any $\varepsilon > 0$, there exists a positive integer $N$ such that
\begin{equation}\label{r26}
\sum\limits_{|j_{1}|\geq N}\sup\limits_{x_{1}\in[0,1]}\sum\limits_{j_{2}\in\mathbf{Z}^{d}} \sup\limits_{x_{2}\in[0,1]^{d}}\Big|\phi(x_{1}+j_{1},x_{2}+j_{2})\Big|\leq \frac{\varepsilon}{4},
\end{equation}
\begin{equation}\label{rr26}
\sum\limits_{|j_{1}|< N}\sup\limits_{x_{1}\in[0,1]}\sum\limits_{|j_{2}|\geq N}\sup\limits_{x_{2}\in[0,1]^{d}}\Big|\phi(x_{1}+j_{1},x_{2}+j_{2})\Big|\leq \frac{\varepsilon}{4}.
\end{equation}
Moreover, we  can obtain
\begin{eqnarray*}
\|\phi_{1}-\phi\|_{\mathcal{W}(L^{1,1})}=&&\sum\limits_{j_{1}\in\mathbf{Z}}\sup\limits_{x_{1}\in[0,1]}\sum\limits_{j_{2}\in\mathbf{Z}^{d}}\sup\limits_{x_{2}\in[0,1]^{d}}\bigg|\phi(x_{1}-j_{1},x_{2}-j_{2})\chi_{O_{N_{0}}}(x_{1}-j_{1},x_{2}-j_{2})\\
&&+\bigg(\sum\limits_{|k_{1}|\geq N}\sum\limits_{|k_{2}|\geq N}+\sum\limits_{|k_{1}|\geq N}\sum\limits_{|k_{2}|< N}+\sum\limits_{|k_{1}|< N}\sum\limits_{|k_{2}|\geq N}\bigg)\phi(x_{1}+k_{1}-j_{1},x_{2}+k_{2}-j_{2})\\
&&\times\chi_{[0,1]\times[0,1]^{d}}(x_{1}-j_{1},x_{2}-j_{2})-\phi(x_{1}-j_{1},x_{2}-j_{2})\bigg|\\
\leq&& F_{1}+F_{2}.
\end{eqnarray*}
Furthermore, we estimate $F_1$ and $F_2$, respectively. In fact, one has
\begin{eqnarray}\label{r29}
F_{1}=&&\sum\limits_{|j_{1}|\geq N}\sup\limits_{x_{1}\in[0,1]}\sum\limits_{j_{2}\in\mathbf{Z}^{d}}\sup\limits_{x_{2}\in[0,1]^{d}}\bigg|\phi(x_{1}-j_{1},x_{2}-j_{2})\chi_{O_{N_{0}}}(x_{1}-j_{1},x_{2}-j_{2})\nonumber\\
&&+\bigg(\sum\limits_{|k_{1}|\geq N}\sum\limits_{|k_{2}|\geq N}+\sum\limits_{|k_{1}|\geq N}\sum\limits_{|k_{2}|< N}+\sum\limits_{|k_{1}|< N}\sum\limits_{|k_{2}|\geq N}\bigg)\phi(x_{1}+k_{1}-j_{1},x_{2}+k_{2}-j_{2})\nonumber\\
&&\times\chi_{[0,1]\times[0,1]^{d}}(x_{1}-j_{1},x_{2}-j_{2})-\phi(x_{1}-j_{1},x_{2}-j_{2})\bigg|\nonumber\\
=&&\sum\limits_{|j_{1}|\geq N} \sup\limits_{x_{1}\in[0,1]}\sum\limits_{j_{2}\in\mathbf{Z}^{d}}\sup\limits_{x_{2}\in[0,1]^{d}}\Big|\phi(x_{1}-j_{1},x_{2}-j_{2})\Big|\nonumber\\
\leq&& \frac{\varepsilon}{4},
\end{eqnarray}
\begin{eqnarray*}
F_{2}=&&\sum\limits_{|j_{1}|< N}\sup\limits_{x_{1}\in[0,1]}\sum\limits_{j_{2}\in\mathbf{Z}^{d}}\sup\limits_{x_{2}\in[0,1]^{d}}\bigg|\phi(x_{1}-j_{1},x_{1}-j_{2})\chi_{O_{N_{0}}}(x_{1}-j_{1},x_{2}-j_{2})\\
&&+\bigg(\sum\limits_{|k_{1}|\geq N}\sum\limits_{|k_{2}|\geq N}+\sum\limits_{|k_{1}|\geq N}\sum\limits_{|k_{2}|< N}+\sum\limits_{|k_{1}|< N}\sum\limits_{|k_{2}|\geq N}\bigg)\phi(x_{1}+k_{1}-j_{1},x_{2}+k_{2}-j_{2})\\
&&\times\chi_{[0,1]\times[0,1]^{d}}(x_{1}-j_{1},x_{2}-j_{2})-\phi(x_{1}-j_{1},x_{2}-j_{2})\bigg|\\
\leq && F_{21}+F_{22},
\end{eqnarray*}
where
\begin{eqnarray}\label{r30}
F_{21}=&&\sum\limits_{|j_{1}|< N}\sup\limits_{x_{1}\in[0,1]}\sum\limits_{|j_{2}|<N}\sup\limits_{x_{2}\in[0,1]^{d}}\bigg|\phi(x_{1}-j_{1},x_{2}-j_{2})\chi_{O_{N_{0}}}(x_{1}-j_{1},x_{2}-j_{2})\nonumber\\
&&+\bigg(\sum\limits_{|k_{1}|\geq N}\sum\limits_{|k_{2}|\geq N}+\sum\limits_{|k_{1}|\geq N}\sum\limits_{|k_{2}|< N}+\sum\limits_{|k_{1}|< N}\sum\limits_{|k_{2}|\geq N}\bigg)\phi(x_{1}+k_{1}-j_{1},x_{2}+k_{2}-j_{2})\nonumber\\
&&\times\chi_{[0,1]\times[0,1]^{d}}(x_{1}-j_{1},x_{2}-j_{2})-\phi(x_{1}-j_{1},x_{2}-j_{2})\bigg|\nonumber\\
=&&\sup\limits_{x_{1}\in[0,1]}\sup\limits_{x_{2}\in[0,1]^{d}}\bigg|
\bigg(\sum\limits_{|k_{1}|\geq N}\sum\limits_{|k_{2}|\geq N}+\sum\limits_{|k_{1}|\geq N}\sum\limits_{|k_{2}|< N}+\sum\limits_{|k_{1}|< N}\sum\limits_{|k_{2}|\geq N}\bigg)\phi(x_{1}+k_{1},x_{2}+k_{2})\bigg|\nonumber\\
\leq&&\sum\limits_{|k_{1}|\geq N}\sup\limits_{x_{1}\in[0,1]}
\sum\limits_{k_{2}\in\mathbf{Z}^{d}}\sup\limits_{x_{2}\in[0,1]^{d}}\Big|\phi(x_{1}+k_{1},x_{2}+k_{2})\Big|+\nonumber\\
&&\sum\limits_{|k_{1}|< N}\sup\limits_{x_{1}\in[0,1]}\sum\limits_{|k_{2}|\geq N}\sup\limits_{x_{2}\in[0,1]^{d}}\Big|\phi(x_{1}+k_{1},x_{2}+k_{2})\Big|\nonumber\\
\leq&&\frac{\varepsilon}{4}+\frac{\varepsilon}{4}=\frac{\varepsilon}{2}
\end{eqnarray}
and
\begin{eqnarray*}
F_{22}=&&\sum\limits_{|j_{1}|< N} \sup\limits_{x_{1}\in[0,1]}\sum\limits_{|j_{2}|\geq N}\sup\limits_{x_{2}\in[0,1]^{d}}\bigg|\phi(x_{1}-j_{1},x_{2}-j_{2})\chi_{O_{N_{0}}}(x_{1}-j_{1},x_{2}-j_{2})\\
&&+\bigg(\sum\limits_{|k_{1}|\geq N}\sum\limits_{|k_{2}|\geq N}+\sum\limits_{|k_{1}|\geq N}\sum\limits_{|k_{2}|< N}+\sum\limits_{|k_{1}|< N}\sum\limits_{|k_{2}|\geq N}\bigg)\phi(x_{1}+k_{1}-j_{1},x_{2}+k_{2}-j_{2})\\
&&\times\chi_{[0,1]\times[0,1]^{d}}(x_{1}-j_{1},x_{2}-j_{2})-\phi(x_{1}-j_{1},x_{2}-j_{2})\bigg|\\
=&&\sum\limits_{|j_{1}|< N}\sup\limits_{x_{1}\in[0,1]}\sum\limits_{|j_{2}|\geq N}\sup\limits_{x_{2}\in[0,1]^{d}}\Big|\phi(x_{1}-j_{1},x_{2}-j_{2})\Big|\\
\leq&&\frac{\varepsilon}{4}.
\end{eqnarray*}
This together with $(\ref{r29})$ and $(\ref{r30})$ obtains
\begin{equation*}
\|\phi_{1}-\phi\|_{\mathcal{W}(L^{1,1})}\leq \varepsilon.
\end{equation*}
Finally, using $(\ref{S11}),\;(\ref{r2.23})$ and $(\ref{r2.25})$ and a similar method as the first case, we can obtain 
\begin{eqnarray*}
&&2^{-n(d+1)}
\Big\|\sum\limits_{j_{1}\in\mathbf{Z}}\sum\limits_{j_{2}\in\mathbf{Z}^{d}}h(2^{-n}j_{1},2^{-n}j_{2})\phi_{1}(x_{1}-j_{1}+k_{1},x_{2}-j_{2}+k_{2})\Big\|_{\mathcal{L}^{\infty,\infty}}\\
\leq&& 2^{-n(d+1)}
\sup\limits_{x_{1}\in[0,1]}\sum\limits_{k_{1}\in\mathbf{Z}}\sup\limits_{x_{2}\in[0,1]^{d}}\sum\limits_{k_{2}\in\mathbf{Z}^{d}}\sum\limits_{j_{1}\in\mathbf{Z}}\sum\limits_{j_{2}\in\mathbf{Z}^{d}}\Big|h(2^{-n}j_{1},2^{-n}j_{2})-h(2^{-n}k_{1},2^{-n}k_{2})\Big|\\
&&\Big|\phi_{1}(x_{1}-j_{1}+k_{1},x_{2}-j_{2}+k_{2})\Big|\\
\leq&& C_{1}(N)C2^{-n(d+1)}2^{-2n}\sup\limits_{x_{1}\in[0,1]}\sum\limits_{k_{1}\in\mathbf{Z}}\sup\limits_{x_{2}\in[0,1]^{d}}\sum\limits_{k_{2}\in\mathbf{Z}^{d}}\sum\limits_{|j_{1}-k_{1}|<N}\sum\limits_{|j_{2}-k_{2}|<N}\Big(1+min(|2^{-n}j_{1}|,|2^{-n}k_{1}|)\Big)^{-1-\varepsilon_{1}}\\
&&\Big(1+min(|2^{-n}j_{2}|,|2^{-n}k_{2}|)\Big)^{-d-\varepsilon_{2}}\Big|\phi_{1}(x_{1}-j_{1}+k_{1},x_{2}-j_{2}+k_{2})\Big|\\
\leq&&
C_{1}(N)C2^{-n(d+1)}2^{-2n}\sup\limits_{x_{1}\in[0,1]}\sum\limits_{k_{1}\in\mathbf{Z}}\Big(1+2^{-n}|k_{1}|\Big)^{-1-\varepsilon_{1}}\sup\limits_{x_{2}\in[0,1]^{d}}\sum\limits_{k_{2}\in\mathbf{Z}^{d}}\Big(1+2^{-n}|k_{2}|\Big)^{-d-\varepsilon_{2}}\\
&&\sum\limits_{|j_{1}-k_{1}|<N}\sum\limits_{|j_{2}-k_{2}|<N}\Big|\phi_{1}(x_{1}-j_{1}+k_{1},x_{2}-j_{2}+k_{2})\Big|\\
\leq&&2^{n(1+\varepsilon_{1})}2^{n(d+\varepsilon_{2})}2^{-n(d+1)}2^{-2n}CC(N)\|\phi_{1}\|_{\mathcal{W}(L^{1,1})}\\
\leq&&2^{n(\varepsilon_{1}-1)}2^{n(\varepsilon_{2}-1)}CC(N)\Big(\|\phi\|_{\mathcal{W}(L^{1,1})}+\varepsilon\Big).
\end{eqnarray*}
This together with the fact that $\mathcal{L}^{\infty,\infty}\subset\mathcal{L}^{p,\infty}$ and $\mathcal{L}^{\infty,\infty}\subset\mathcal{L}^{\infty,q}$ proves  $(\ref{rr2.15})$ for the cases $p=\infty$ or $q=\infty$.

\section{Main results}
\ \ \ \ In this section, we will give the main result of this paper and its proofs.
\begin{Tm}\label{M3.1}
Let $\Phi=(\phi_{1},\cdots,\phi_{r})^{\mathrm{T}}\in \mathcal{L}^{\infty,\infty}$ if $1 < p,q < \infty$ and $\Phi \in \mathcal{W}(L^{1,1})$ if $p=1,\infty$ or $q=1,\infty$. Then the following statements are equivalent to each other.
\par $(i)\;\;V_{p,q}(\Phi)$ is closed in $L^{p,q}(\mathbf{R}\times \mathbf{R}^{d})$.
\par $(ii)\;\left\{\phi_{i}(\cdot-j_{1},\cdot-j_{2}):(j_{1},j_{2})\in \mathbf{Z}\times\mathbf{Z}^{d},1\leq i \leq r\right\}$ is a $(p,q)$-frame for $V_{p,q}(\Phi)$, i.e., there exists a positive constant A (depending on $p,\;q$ and $\Phi$) such that
\begin{equation}\label{M1}
A^{-1}\|f\|_{L^{p,q}}\leq \sum\limits_{i=1}^{r}\left\|\bigg(\int_{\mathbf{R}}\int _{\mathbf{R}^{d}}f(x_{1},x_{2})\overline{\phi_{i}(x_{1}-j_{1},x_{2}-j_{2})}dx_{1}dx_{2}\bigg)_{(j_{1},j_{2})\in \mathbf{Z}\times \mathbf{Z}^{d}}\right\|_{\ell^{p,q}}\leq A\|f\|_{L^{p,q}}
\end{equation}
holds for any $f\in V_{p,q}(\Phi)$.
\par $(iii)$  There exists a positive constant C such that\\
\begin{equation*}
C^{-1}[\hat{\Phi},\hat{\Phi}](\xi,\tilde{\xi})
\leq [\hat{\Phi},\hat{\Phi}](\xi,\tilde{\xi})\overline{[\hat{\Phi},\hat{\Phi}](\xi,\tilde{\xi})^{\mathrm{T}}}\leq C[\hat{\Phi},\hat{\Phi}](\xi,\tilde{\xi}),\;\;\;\; (\xi,\tilde{\xi})\in [-\pi,\pi]\times[-\pi,\pi]^{d}.
\end{equation*}
\par $(iv)$  There exists a positive constant B (depending on p,q and $\Phi$) such that
\begin{equation}\label{r3.2}
B^{-1}\|f\|_{L^{p,q}}\leq \inf\limits_{f=\sum\limits_{i=1}^{r}\phi_{i}\ast^{'}D_{i}}\sum\limits_{i=1}^{r}\|D_{i}\|_{\ell^{p,q}}\leq B\|f\|_{L^{p,q}}\;\;\;for\;f\in V_{p,q}(\Phi).
\end{equation}
\par $(v)$ There exists $\Psi=(\psi_{1},\cdots,\psi_{r})^{\mathrm{T}}\in \mathcal{L}^{\infty,\infty}$ if $1 < p,q < \infty$ and $\Psi\in \mathcal{W}(L^{1,1})$ if $p=1,\infty$ or $q=1,\infty$ such that for any $f\in V_{p,q}(\Phi)$,
\begin{equation}\label{r3.3}
\begin{split}
f &= \sum\limits_{i=1}^{r}\sum\limits_{j_{1}\in \mathbf{Z}}\sum\limits_{j_{2}\in \mathbf{Z}^{d}}\langle f,\psi_{i}(\cdot-j_{1},\cdot-j_{2})\rangle \phi_{i}(\cdot-j_{1},\cdot-j_{2})\\
&= \sum\limits_{i=1}^{r}\sum\limits_{j_{1}\in \mathbf{Z}}\sum\limits_{j_{2}\in \mathbf{Z}^{d}}\langle f,\phi_{i}(\cdot-j_{1},\cdot-j_{2})\rangle \psi_{i}(\cdot-j_{1},\cdot-j_{2}).
 \end{split}
\end{equation}
\end{Tm}
\subsection{Proof of $ (v)\Longrightarrow (iv)$}
Let $f= \sum\limits_{i=1}^{r}\sum\limits_{j_{1}\in \mathbf{Z}}\sum\limits_{j_{2}\in \mathbf{Z}^{d}}\langle f,\psi_{i}(\cdot-j_{1},\cdot-j_{2})\rangle \phi_{i}(\cdot-j_{1},\cdot-j_{2})$. Then using Lemma $\ref{L2.5}$, we obtain
\begin{eqnarray}\label{r3.4}
\inf\limits_{f=\sum\limits_{i=1}^{r}\phi_{i}\ast^{'}D_{i}}\sum\limits_{i=1}^{r}\|D_{i}\|_{\ell^{p,q}}&&\leq \sum\limits_{i=1}^{r}\Big\|\Big\{<f,\psi_{i}(\cdot-j_{1},\cdot-j_{2})>\Big\}_{(j_{1},j_{2})\in \mathbf{Z}\times\mathbf{Z}^{d}}\Big\|_{\ell^{p,q}}\nonumber\\
&&\leq \bigg(\sum\limits_{i=1}^{r}\|\psi_{i}\|_{\mathcal{L}^{\infty,\infty}}\bigg)\|f\|_{L^{p,q}}.
\end{eqnarray}
\par If $f=\sum\limits_{i=1}^{r}\phi_{i}\ast^{'}D_{i}\in V_{p,q}(\Phi)$, then by Lemma $\ref{Lm 2.1}$, we obtain
\begin{equation*}
\|f\|_{L^{p,q}}=\bigg\|\sum\limits_{i=1}^{r}\phi_{i}\ast^{'}D_{i}\bigg\|_{L^{p,q}}\leq\sum_{i=1}^{r}\Big\|\phi_{i}\ast^{'}D_{i}\Big\|_{L^{p,q}} \leq\bigg(\sum_{i=1}^{r}\|\phi_{i}\|_{\mathcal{L}^{\infty,\infty}}\bigg)\bigg(\sum_{i=1}^{r}\|D_{i}\|_{\ell^{p,q}}\bigg).
\end{equation*}
This together with $(\ref{r3.4})$ proves the inequality $(\ref{r3.2})$ for
$B=\max\bigg\{\sum\limits_{i=1}^{r}\|\psi_{i}\|_{\mathcal{L}^{\infty,\infty}},\sum\limits_{i=1}^{r}\|\phi_{i}\|_{\mathcal{L}^{\infty,\infty}}\bigg\}$.
\subsection{Proof of $(iv)\Longrightarrow (i)$}
\begin{Lm}\cite{A1}\label{L3.2}
Let $(X,\|\cdot\|_{X})$ and $(Y,\|\cdot\|_{Y})$ be two Banach spaces, $T$ be a bounded linear operator from $X$ to $Y$. If there exists a positive constant $C$ such that \\
\begin{equation*}
C^{-1}\|y\|_{Y}\leq \inf\limits_{y=Tx}\|x\|_{X}\leq C\|y\|_{Y}\;\;\;for\;y\in Ran(T),
\end{equation*}
then $Ran(T)$ is closed.\\
Proof of $(iv)\Longrightarrow (i)$. Take $X=(\ell^{p,q})^{r}$ and $Y=L^{p,q}(\mathbf{R}\times\mathbf{R}^{d})$. Define $T:\;X\rightarrow Y$ by\\
\begin{equation*}
(Tc)(x,y)=\sum\limits_{i=1}^{r}\sum\limits_{j_{1}\in\mathbf{Z}}\sum\limits_{j_{2}\in\mathbf{Z}^{d}}c_{i}(j_{1},j_{2})\phi_{i}(x-j_{1},y-j_{2})
\end{equation*}
for $c=(c_{1},\cdots,c_{r})^{\mathrm{T}}\in(\ell^{p,q})^{r}.$ Then $Ran(T)=V_{p,q}(\Phi)$.
Finally, it follows from Lemma $\ref{L3.2}$ and $(\ref{r3.2})$ that $(i)$ holds.
\end{Lm}
\subsection{Proof of $(i)\Rightarrow(iii)$}
Let $k_{0}=\min\limits_{(\xi,\tilde{\xi})\in \mathbf{R}\times\mathbf{R}^{d}}rank\Big(\hat{\Phi}(\xi+2k_{1}\pi,\tilde{\xi}+2k_{2}\pi)\Big)_{(k_{1},k_{2})\in \mathbf{Z}\times \mathbf{Z}^{d}}$ and  $\Omega_{k_{0}}=\Big\{(\xi,\tilde{\xi}):\;rank\Big(\hat{\Phi}(\xi+2k_{1}\pi,\tilde{\xi}+2k_{2}\pi)\Big)_{(k_{1},k_{2})\in \mathbf{Z}\times\mathbf{Z}^{d}}>k_{0}\Big\}$.
Then $\Omega_{k_{0}}\neq \mathbf{R}^{d+1}$. By Lemma $\ref{L7}$, it suffices to prove that $\Omega_{k_{0}}=\varnothing$.
On the contrary, suppose that $\Omega_{k_{0}}\neq \varnothing$. Then it suffices to construct a function $G$ in the $L^{p,q}$-closure of $V_{p,q}(\Phi)$ such that $G$ cannot be written as $\sum\limits_{i=1}^{r}\phi_{i}\ast^{'}D_{i}$ for some $D_{i}\in \ell^{p,q},\;1\leq i \leq r .$
Note that $\Omega_{k_{0}}$ is an open set. Then the boundary $\partial\Omega_{k_{0}}$ of $\Omega_{k_{0}}$ is nonempty. For any $(\xi_{0},\tilde{\xi}_{0})\in \partial\Omega_{k_{0}},\;rank\Big(\hat{\Phi}(\xi_{0}+2k_{1}\pi,\tilde{\xi}_{0}+2k_{2}\pi)\Big)_{(k_{1},k_{2})\in \mathbf{Z}\times\mathbf{Z}^{d}}=k_{0}$ and\\
\begin{equation*}
\max_{(\xi,\tilde{\xi})\in B((\xi_{0},\tilde{\xi}_{0}),\delta)}rank\bigg(\hat{\Phi}(\xi+2k_{1}\pi,\tilde{\xi}+2k_{2}\pi)\bigg)_{(k_{1},k_{2})\in \mathbf{Z}\times \mathbf{Z}^{d}}>k_{0}\;\;\;for\;\;\delta>0.
\end{equation*}
Similar to the proof of Lemma $\ref{L8}$, there exists a $r\times r$ $2\pi$-periodic nonsingular matrix $P_{(\xi_{0},\tilde{\xi}_{0})}(\xi,\tilde{\xi})$ and $(K_{\xi_{0}},K_{\tilde{\xi}_{0}})\subset\mathbf{Z}\times\mathbf{Z}^{d}$ with cardinality $\#(K_{\xi_{0}},K_{\tilde{\xi}_{0}})=k_{0}$ such that $\Psi_{(\xi_{0},\tilde{\xi}_{0})}$ defined by \\
\begin{equation*}
\hat{\Psi}_{(\xi_{0},\tilde{\xi}_{0})}(\xi,\tilde{\xi})=
\begin{pmatrix}
\hat{\Psi}_{1,(\xi_{0},\tilde{\xi}_{0})}(\xi,\tilde{\xi})\\
\hat{\Psi}_{2,(\xi_{0},\tilde{\xi}_{0})}(\xi,\tilde{\xi})
\end{pmatrix}
=P_{(\xi_{0},\tilde{\xi}_{0})}(\xi,\tilde{\xi})\hat{\Phi}(\xi,\tilde{\xi})
\end{equation*}
satisfies
\begin{equation}\label{M2}
rank\bigg(\hat{\Psi}_{1,(\xi_{0},\tilde{\xi}_{0})}(\xi+2k_{1}\pi,\tilde{\xi}+2k_{2}\pi)\bigg)_{(k_{1},k_{2})\in (K_{\xi_{0}}, K_{\tilde{\xi}_{0}})}=k_{0}\;for\;(\xi,\tilde{\xi})\in B\Big((\xi_{0},\tilde{\xi}_{0}),2\delta_{(\xi_{0},\tilde{\xi}_{0})}\Big),\;\;\;
\end{equation}
\begin{equation}\label{M3}
\hat{\Psi}_{2,(\xi_{0},\tilde{\xi}_{0})}(\xi_{0}+2k_{1}\pi,\tilde{\xi}_{0}+2k_{2}\pi)=0\;for\; (k_{1},k_{2})\in\mathbf{Z}\times\mathbf{Z}^{d},
\end{equation}
\begin{equation}\label{M4}
\hat{\Psi}_{2,(\xi_{0},\tilde{\xi}_{0})}(\xi+2k_{1}\pi,\tilde{\xi}+2k_{2}\pi)=0\;for\;(k_{1}, k_{2})\in (K_{\xi_{0}},K_{\tilde{\xi}_{0}})\;and\;(\xi,\tilde{\xi})\in B\Big((\xi_{0},\tilde{\xi}_{0}),2\delta_{(\xi_{0},\tilde{\xi}_{0})}\Big)\;
\end{equation}
and
\begin{equation}\label{M5}
\hat{\Psi}_{2,(\xi_{0},\tilde{\xi}_{0})}(\xi,\tilde{\xi})\not\equiv 0\;\;on \;\;B\Big((\xi_{0},\tilde{\xi}_{0}),\delta\Big)+2\pi(\mathbf{Z}\times\mathbf{Z}^{d})
\end{equation}
for all $0<\delta< 2\delta_{(\xi_{0},\tilde{\xi}_{0})}$. Since the Fourier coefficients of $P_{(\xi_{0},\tilde{\xi}_{0})}(\xi,\tilde{\xi})$ belong to $\ell^{1,1}$ and $\Phi\in\mathcal{L}^{\infty,\infty}$ (or $\Phi\in \mathcal{W}(L^{1,1})$), we have $\Psi_{(\xi_{0},\tilde{\xi}_{0})}\in \mathcal{L}^{\infty,\infty}$ (or $\Psi_{(\xi_{0},\tilde{\xi}_{0})}\in \mathcal{W}(L^{1,1})$). This together with $(\ref{M3})$ and the Poisson's summation formula leads to
\begin{eqnarray}\label{M6}
&&\sum\limits_{j_{1}\in\mathbf{Z}}\sum\limits_{j_{2}\in\mathbf{Z}^{d}} e^{-i(\xi_{0},\tilde{\xi}_{0})\cdot(\cdot+j_{1},\cdot+j_{2})} \Psi_{2,(\xi_{0},\tilde{\xi}_{0})}(\cdot+j_{1},\cdot+j_{2})\nonumber\\.
&&=\sum\limits_{j_{1}\in\mathbf{Z}}\sum\limits_{j_{2}\in\mathbf{Z}^{d}}\hat{\Psi}_{2,(\xi_{0},\tilde{\xi}_{0})}(\xi_{0}+2\pi j_{1},\tilde{\xi}_{0}+2\pi j_{2})e^{-i2\pi(j_{1},j_{2})\cdot(\cdot,\cdot)}\nonumber\\
&&=0.
\end{eqnarray}
Let $H$ be a nonnegative $C^{\infty}$ function satisfying $(\ref{S7})$ and
\begin{equation*}
H_{n,(\xi_{0},\tilde{\xi}_{0})}(\xi,\tilde{\xi})=\sum\limits_{k_{1}\in\mathbf{Z}}\sum\limits_{k_{2}\in\mathbf{Z}^{d}}H\Big(2^{n}(\xi+2k_{1}\pi-\xi_{0}),2^{n}(\tilde{\xi}+2k_{2}\pi-\tilde{\xi}_{0})\Big),\;\;n\geq2.
\end{equation*}
Define an operator $T_{n,(\xi_{0},\tilde{\xi}_{0})}$ from $(\ell^{p,q})^{r-k_{0}}$ to $L^{p,q}$ by
\begin{equation*}
T_{n,(\xi_{0},\tilde{\xi}_{0})}:D=(D_{1},\cdots,D_{r-k_{0}})^{\mathrm{T}}\rightarrow \sum\limits_{i=1}^{r-k_{0}}\Big(\psi_{2,(\xi_{0},\tilde{\xi}_{0}),i}\ast^{'}h_{n,(\xi_{0},\tilde{\xi}_{0})}\Big)\ast^{'}D_{i},
\end{equation*}
where $h_{n,(\xi_{0},\tilde{\xi}_{0})}$ denotes the sequence of the Fourier coefficients of the $2\pi$-periodic function $H_{n,(\xi_{0},\tilde{\xi}_{0})}(\xi,\tilde{\xi})$ and $\Psi_{2,(\xi_{0},\tilde{\xi}_{0})}=(\psi_{2,(\xi_{0},\tilde{\xi}_{0}),1},\cdots,\psi_{2,(\xi_{0},\tilde{\xi}_{0}),r-k_{0}})^{\mathrm{T}}$. By $(\ref{M6})$, Lemma $\ref{Lm 2.1}$ and Lemma $\ref{L9}$, we obtain
\begin{eqnarray}\label{r9}
\|T_{n,(\xi_{0},\tilde{\xi}_{0})}D\|_{L^{p,q}}&&\leq\sum\limits_{i=1}^{r-k_{0}}\bigg\|\Big(\psi_{2,(\xi_{0},\tilde{\xi}_{0}),i}\ast'h_{n,(\xi_{0},\tilde{\xi}_{0})}\Big)\ast'D_{i}\bigg\|_{L^{p,q}}\nonumber\\
&&\leq\sum\limits_{i=1}^{r-k_{0}}\|D_{i}\|_{\ell^{p,q}}\|\psi_{2,(\xi_{0},\tilde{\xi}_{0}),i}\ast'h_{n,(\xi_{0},\tilde{\xi}_{0})}\|_{\mathcal{L}^{p,q}}\nonumber\\
&&=\sum\limits_{i=1}^{r-k_{0}}\|D_{i}\|_{\ell^{p,q}}\Big\|\sum\limits_{j_{1}\in\mathbf{Z}}\sum\limits_{j_{2}\in\mathbf{Z}^{d}}h_{n,(\xi_{0},\tilde{\xi}_{0})}(j_{1},j_{2})\psi_{2,(\xi_{0},\tilde{\xi}_{0}),i}(\cdot-j_{1},\cdot-j_{2})\Big\|_{\mathcal{L}^{p,q}}.\;\;\;\;\;\;\;\;\;\;
\end{eqnarray}
Moreover, we can obtain
\begin{eqnarray*}
h_{n,(\xi_{0},\tilde{\xi}_{0})}(j_{1},j_{2})&&=\int_{-\pi}^{\pi}\int_{[-\pi,\pi]^{d}}H_{n,(\xi_{0},\tilde{\xi}_{0})}(\xi,\tilde{\xi})e^{-i(j_{1},j_{2})\cdot(\xi,\tilde{\xi})}d\xi d\tilde{\xi}\\
&&=\sum\limits_{k_{1}\in\mathbf{Z}}\sum\limits_{k_{2}\in\mathbf{Z}^{d}}\int_{-\pi}^{\pi}\int_{[-\pi,\pi]^{d}}H\Big(2^{n}(\xi+2k_{1}\pi-\xi_{0}),2^{n}(\tilde{\xi}+2k_{2}\pi-\tilde{\xi}_{0})\Big)e^{-i(j_{1},j_{2})\cdot(\xi,\tilde{\xi})}d\xi d\tilde{\xi}\\
&&=e^{-i(j_{1},j_{2})\cdot(\xi_{0},\tilde{\xi}_{0})}\int_{\mathbf{R}}\int_{\mathbf{R}^{d}}H(2^{n}\xi,2^{n}\tilde{\xi})e^{-i(j_{1},j_{2})\cdot(\xi,\tilde{\xi})}d\xi d\tilde{\xi}\\
&&=2^{-n(d+1)}e^{-i(j_{1},j_{2})\cdot(\xi_{0},\tilde{\xi}_{0})}\int_{\mathbf{R}}\int_{\mathbf{R}^{d}}H(\xi,\tilde{\xi})e^{-i(j_{1},j_{2})\cdot(2^{-n}\xi,2^{-n}\tilde{\xi})}d\xi d\tilde{\xi}\\
&&=e^{-i(j_{1},j_{2})\cdot(\xi_{0},\tilde{\xi}_{0})}2^{-n(d+1)}\hat{H}(2^{-n}j_{1},2^{-n}j_{2}).
\end{eqnarray*}
This together with $(\ref{M6})$ and $(\ref{r9})$ obtains
\begin{equation*}
\|T_{n,(\xi_{0},\tilde{\xi}_{0})}D\|_{L^{p,q}}\leq\|D\|_{\ell^{p,q}}2^{-n(d+1)}\Big\|\sum\limits_{j_{1}\in\mathbf{Z}}\sum\limits_{j_{2}\in\mathbf{Z}^{d}} \hat{H}(2^{-n}j_{1},2^{-n}j_{2})e^{-i(\xi_{0},\tilde{\xi}_{0})\cdot(\cdot+j_{1},\cdot+j_{2})}\Psi_{2,(\xi_{0},\tilde{\xi}_{0})}(\cdot+j_{1},\cdot+j_{2})\Big\|_{\mathcal{L}^{p,q}}.\\
\end{equation*}
Then it follows from Lemma $\ref{L9}$ that
\begin{equation}\label{r10}
\lim\limits_{n\rightarrow\infty}\|T_{n,(\xi_{0},\tilde{\xi}_{0})}\|=0.
\end{equation}
 Let $\tilde{H}(x_{1},x_{2})=H(2x_{1},2x_{2})-H(8x_{1},8x_{2})$ and define
\begin{equation*}
\tilde{H}_{n,(\xi_{0},\tilde{\xi}_{0})}(\xi,\tilde{\xi})=\sum\limits_{k_{1}\in\mathbf{Z}}\sum\limits_{k_{2}\in\mathbf{Z}^{d}}\tilde{H}\Big(2^{n}(\xi+2k_{1}\pi-\xi_{0}),2^{n}(\tilde{\xi}+2k_{2}\pi-\tilde{\xi}_{0})\Big),\;\;n\geq2.
\end{equation*}
Then for $n\geq 4$, we have
\begin{equation}\label{M8}
\tilde{H}_{n,(\xi_{0},\tilde{\xi}_{0})}(\xi,\tilde{\xi})H_{n,(\xi_{0},\tilde{\xi}_{0})}(\xi,\tilde{\xi})=\tilde{H}_{n,(\xi_{0},\tilde{\xi}_{0})}(\xi,\tilde{\xi}).
\end{equation}
Define an operator $\tilde{T}_{n,(\xi_{0},\tilde{\xi}_{0})}$ from $(\ell^{p,q})^{r-k_{0}}$ to $L^{p,q}$ by
\begin{equation}
\tilde{T}_{n,(\xi_{0},\tilde{\xi}_{0})}:D=(D_{1},\cdots,D_{r-k_{0}})^{\mathrm{T}}\rightarrow \sum\limits_{i=1}^{r-k_{0}}\Big(\psi_{2,(\xi_{0},\tilde{\xi}_{0}),i}\ast^{'}\tilde{h}_{n,(\xi_{0},\tilde{\xi}_{0})}\Big)\ast^{'}D_{i},
\end{equation}
where $\tilde{h}_{n,(\xi_{0},\tilde{\xi}_{0})}$ denotes the sequence of the Fourier coefficients of the $2\pi$-periodic function $\tilde{H}_{n,(\xi_{0},\tilde{\xi}_{0})}(\xi,\tilde{\xi})$.
Note that
\begin{equation*}
\tilde{H}_{n,(\xi_{0},\tilde{\xi}_{0})}(\xi,\tilde{\xi})=H_{n+1,(\xi_{0},\tilde{\xi}_{0})}(\xi,\tilde{\xi})-H_{n+3,(\xi_{0},\tilde{\xi}_{0})}(\xi,\tilde{\xi}).
\end{equation*}
Then it follows from $(\ref{r10})$ that
\begin{eqnarray}
\|\tilde{T}_{n,(\xi_{0},\tilde{\xi}_{0})}\|&&=\|T_{n+1,(\xi_{0},\tilde{\xi}_{0})}-T_{n+3,(\xi_{0},\tilde{\xi}_{0})}\|\nonumber\\
&&\leq\|T_{n+1,(\xi_{0},\tilde{\xi}_{0})}\|+\|T_{n+3,(\xi_{0},\tilde{\xi}_{0})}\|\nonumber\\
&&\rightarrow 0,\;as\;n\rightarrow\infty.
\end{eqnarray}
By $(\ref{M5})$ and $(\ref{r10})$, there exists a subsequence $n_{\ell},\;\ell\geq 1$ such that $n_{\ell+1}\geq n_{\ell}+8$,
\begin{equation}\label{M9}
\|\tilde{T}_{n_{\ell},(\xi_{0},\tilde{\xi}_{0})}\|\neq 0\;\;and \;\;\sum\limits_{k=0}^{5}\|T_{n_{\ell}+k,(\xi_{0},\tilde{\xi}_{0})}\|\leq 2^{-\ell}.
\end{equation}
Let $D_{n_{\ell}}\in(\ell^{p,q})^{r-k_{0}}$ be chosen such that
\begin{equation}\label{r3.16}
\|D_{n_{\ell}}\|_{\ell^{p,q}}=1\;\;and \;\;\|\tilde{T}_{n_{\ell},(\xi_{0},\tilde{\xi}_{0})}D_{n_{\ell}}\|_{L^{p,q}}\geq \|\tilde{T}_{n_{\ell},(\xi_{0},\tilde{\xi}_{0})}\|/2.
\end{equation}
For sufficiently large $\ell_{1}$ and $s\geq \ell_{1}$, define $G_{s}$ and $G$ by
\begin{equation*}
\begin{split}
\hat{G}_{s}(\xi,\tilde{\xi})&=\sum\limits_{\ell=\ell_{1}}^{s}\ell\Big(H_{n_{\ell},(\xi_{0},\tilde{\xi}_{0})}(\xi,\tilde{\xi})-H_{n_{\ell}+4,(\xi_{0},\tilde{\xi}_{0})}(\xi,\tilde{\xi})\Big)D_{n_{\ell}}(\xi,\tilde{\xi})^{\mathrm{T}}\hat{\Psi}_{2,(\xi_{0},\tilde{\xi}_{0})}(\xi,\tilde{\xi})\\
&=\sum\limits_{\ell=\ell_{1}}^{s}\ell\Big(H_{n_{\ell},(\xi_{0},\tilde{\xi}_{0})}(\xi,\tilde{\xi})-H_{n_{\ell}+4,(\xi_{0},\tilde{\xi}_{0})}(\xi,\tilde{\xi})\Big)\Big(0,D _{n_{\ell}}(\xi,\tilde{\xi})^{\mathrm{T}}\Big)P_{(\xi_{0},\tilde{\xi}_{0})}(\xi,\tilde{\xi})\hat{\Phi}(\xi,\tilde{\xi})
\end{split}
\end{equation*}
and
\begin{equation*}
\hat{G}(\xi,\tilde{\xi})=\sum\limits_{\ell=\ell_{1}}^{\infty}\ell\Big(H_{n_{\ell},(\xi_{0},\tilde{\xi}_{0})}(\xi,\tilde{\xi})-H_{n_{\ell}+4,(\xi_{0},\tilde{\xi}_{0})}(\xi,\tilde{\xi})\Big)D_{n_{\ell}}(\xi,\tilde{\xi})^{\mathrm{T}}\hat{\Psi}_{2,(\xi_{0},\tilde{\xi}_{0})}(\xi,\tilde{\xi}).
\end{equation*}
 Now it remains to prove that $G$ is in the $L^{p,q}$-closure of $V_{p,q}(\Phi)$ and that $G$ cannot be written as $\sum\limits_{i=1}^{r}\phi_{i}\ast^{'}D_{i}$ for some $D_{i}\in \ell^{p,q},\;1\leq i\leq r.$ From the construction of $G_{s}$ and $G,\;G_{s}\in V_{p,q}(\Phi)$ for all $s\geq \ell_{1}$ and
\begin{eqnarray}\label{r17}
\|G-G_{s}\|_{L^{p,q}}
\leq&& \sum\limits_{\ell=s+1}^{\infty}\ell\bigg\|\mathcal{F}^{-1}\bigg(\Big(H_{n_{\ell},(\xi_{0},\tilde{\xi}_{0})}(\xi,\tilde{\xi})-H_{n_{\ell}+4,(\xi_{0},\tilde{\xi}_{0})}(\xi,\tilde{\xi})\Big)D_{n_{\ell}}(\xi,\tilde{\xi})^{\mathrm{T}}\hat{\Psi}_{2,(\xi_{0},\tilde{\xi}_{0})}(\xi,\tilde{\xi})\bigg)\bigg\|_{L^{p,q}}\nonumber\\
\leq&&\sum\limits_{\ell=s+1}^{\infty}\ell\bigg(\bigg\|\mathcal{F}^{-1}\Big(H_{n_{\ell},(\xi_{0},\tilde{\xi}_{0})}(\xi,\tilde{\xi})D_{n_{\ell}}(\xi,\tilde{\xi})^{\mathrm{T}}\hat{\Psi}_{2,(\xi_{0},\tilde{\xi}_{0})}(\xi,\tilde{\xi})\Big)\bigg\|_{L^{p,q}}\nonumber\\
&&+\bigg\|\mathcal{F}^{-1}\Big(H_{n_{\ell}+4,(\xi_{0},\tilde{\xi}_{0})}(\xi,\tilde{\xi})D_{n_{\ell}}(\xi,\tilde{\xi})^{\mathrm{T}}\hat{\Psi}_{2,(\xi_{0},\tilde{\xi}_{0})}(\xi,\tilde{\xi})\Big)\bigg\|_{L^{p,q}}\bigg).
\end{eqnarray}
Note that
\begin{eqnarray*}
T_{n_{\ell},(\xi_{0},\tilde{\xi}_{0})}D_{n_{\ell}}&&=\sum\limits_{i=1}^{r-k_{0}}\Big(\psi_{2,(\xi_{0},\tilde{\xi}_{0}),i}\ast'h_{n_{\ell},(\xi_{0},\tilde{\xi}_{0})}\Big)\ast'D_{n_{\ell}}^{(i)}\\
&&=\sum\limits_{i=1}^{r-k_{0}}\sum\limits_{k_{1}\in\mathbf{Z}}\sum\limits_{k_{2}\in\mathbf{Z}^{d}}D_{n_{\ell}}^{(i)}(k_{1},k_{2})\Big(\psi_{2,(\xi_{0},\tilde{\xi}_{0}),i}\ast'h_{n_{\ell},(\xi_{0},\tilde{\xi}_{0})}\Big)(x_{1}-k_{1},x_{2}-k_{2}).
\end{eqnarray*}
Then, we obtain
\begin{eqnarray*}
(T_{n_{\ell},(\xi_{0},\tilde{\xi}_{0})}D_{n_{\ell}})^{\wedge}(\xi,\tilde{\xi})&&=\sum\limits_{i=1}^{r-k_{0}}\sum\limits_{k_{1}\in\mathbf{Z}}\sum\limits_{k_{2}\in\mathbf{Z}^{d}}D_{n_{\ell}}^{(i)}(k_{1},k_{2})e^{-i(k_{1},k_{2})\cdot(\xi,\tilde{\xi})}\Big(\psi_{2,(\xi_{0},\tilde{\xi}_{0}),i}\ast'h_{n_{\ell},(\xi_{0},\tilde{\xi}_{0})}\Big)^{\wedge}(\xi,\tilde{\xi})\\
&&=\sum\limits_{i=1}^{r-k_{0}}D_{n_{\ell}}^{(i)}(\xi,\tilde{\xi})\bigg(\sum\limits_{k_{1}\in\mathbf{Z}}\sum\limits_{k_{2}\in\mathbf{Z}^{d}}h_{n_{\ell},(\xi_{0},\tilde{\xi}_{0})}(k_{1},k_{2})\psi_{2,(\xi_{0},\tilde{\xi}_{0}),i}(x_{1}-k_{1},x_{2}-k_{2})\bigg)^{\wedge}(\xi,\tilde{\xi})\\
&&=\sum\limits_{i=1}^{r-k_{0}}D_{n_{\ell}}^{(i)}(\xi,\tilde{\xi})H_{n_{\ell},(\xi_{0},\tilde{\xi}_{0})}(\xi,\tilde{\xi})\hat{\psi}_{2,(\xi_{0},\tilde{\xi}_{0}),i}(\xi,\tilde{\xi})\\
&&=H_{n_{\ell},(\xi_{0},\tilde{\xi}_{0})}(\xi,\tilde{\xi})D_{n_{\ell}}(\xi,\tilde{\xi})^{\mathrm{T}}\hat{\Psi}_{2,(\xi_{0},\tilde{\xi}_{0})}(\xi,\tilde{\xi}).
\end{eqnarray*}
This together with $(\ref{r10})$ and $(\ref{r17})$ obtains
\begin{eqnarray*}
\|G-G_{s}\|_{L^{p,q}}&&\leq \sum\limits_{\ell=s+1}^{\infty}\ell\Big(\|T_{n_{\ell}+4,(\xi_{0},\tilde{\xi}_{0})}\|+\|T_{n_{\ell},(\xi_{0},\tilde{\xi}_{0})}\|\Big)\\
&&\leq\sum\limits_{l=s+1}^{\infty}\frac{\ell}{2^{\ell}}\rightarrow0,\;as\;s\rightarrow\infty.
\end{eqnarray*}
This shows that $G$ is in the $L^{p,q}$-closure of $V_{p,q}(\Phi)$.
\par Finally, we prove that $G\not\in V_{p,q}(\Phi)$. On the contrary, suppose that
\begin{equation}\label{M10}
\hat{G}(\xi,\tilde{\xi})=B(\xi,\tilde{\xi})^{\mathrm{T}}\hat{\Phi}(\xi,\tilde{\xi})
\end{equation}
for some vector valued $2\pi$-periodic distribution $B(\xi,\tilde{\xi})$ with Fourier coefficients belonging to $\ell^{p,q}$. Note that supp$\hat{G}_{s}(\xi,\tilde{\xi})\subset B\Big((\xi_{0},\tilde{\xi}_{0}),2^{-(\ell_{1}-1)}\Big)+2\pi(\mathbf{Z}\times\mathbf{Z}^{d})$ for all $s\geq \ell_{1}$ and so is supp$\hat{G}(\xi,\tilde{\xi})$. Hence we may assume that $B(\xi,\tilde{\xi})$ in $(\ref{M10})$ is supported in $B\Big((\xi_{0},\tilde{\xi}_{0}),\delta_{(\xi_{0},\tilde{\xi}_{0})}\Big)+2\pi(\mathbf{Z}\times\mathbf{Z}^{d})$ when $\ell_{1}$ is chosen large enough. Write
\begin{equation*}
B(\xi,\tilde{\xi})^{\mathrm{T}}\Big(P_{(\xi_{0},\tilde{\xi}_{0})}(\xi,\tilde{\xi})\Big)^{-1}=\bigg(B_{1}(\xi,\tilde{\xi})^{\mathrm{T}},B_{2}(\xi,\tilde{\xi})^{\mathrm{T}}\bigg).
\end{equation*}
Then we can rewrite $(\ref{M10})$ as\\
$B_{1}(\xi,\tilde{\xi})^{\mathrm{T}}\hat{\Psi}_{1,(\xi_{0},\tilde{\xi}_{0})}(\xi,\tilde{\xi})$
\begin{equation}\label{M11}
=\bigg(-B_{2}(\xi,\tilde{\xi})^{\mathrm{T}}+\sum\limits_{\ell=\ell_{1}}^{\infty}\ell\Big(H_{n_{\ell},(\xi_{0},\tilde{\xi}_{0})}(\xi,\tilde{\xi})-H_{n_{\ell}+4,(\xi_{0},\tilde{\xi}_{0})}(\xi,\tilde{\xi})\Big)D_{n_{\ell}}(\xi,\tilde{\xi})^{\mathrm{T}}\bigg)\hat{\Psi}_{2,(\xi_{0},\tilde{\xi}_{0})}(\xi,\tilde{\xi}).
\end{equation}
This together with $(\ref{M2})$ and $(\ref{M4})$ obtains that $B_{1}(\xi,\tilde{\xi})\equiv 0$. Substituting this into $(\ref{M11}),$
\begin{eqnarray}\label{M12}
&&B_{2}(\xi,\tilde{\xi})^{\mathrm{T}}\hat{\Psi}_{2,(\xi_{0},\tilde{\xi}_{0})}(\xi,\tilde{\xi})\nonumber\\
&&=\sum\limits_{\ell=\ell_{1}}^{\infty}\ell\Big(H_{n_{\ell},(\xi_{0},\tilde{\xi}_{0})}(\xi,\tilde{\xi})-H_{n_{\ell}+4,(\xi_{0},\tilde{\xi}_{0})}(\xi,\tilde{\xi})\Big)D_{n_{\ell}}(\xi,\tilde{\xi})^{\mathrm{T}}\hat{\Psi}_{2,(\xi_{0},\tilde{\xi}_{0})}(\xi,\tilde{\xi}).
\end{eqnarray}
By direct computation and using the fact that $n_{\ell+1}\geq n_{\ell}+8$, we obtain
\begin{equation*}
\tilde{H}_{n_{\ell},(\xi_{0},\tilde{\xi}_{0})}(\xi,\tilde{\xi})\Big(H_{n_{\ell^{'}},(\xi_{0},\tilde{\xi}_{0})}(\xi,\tilde{\xi})-H_{n_{\ell^{'}}+4,(\xi_{0},\tilde{\xi}_{0})}(\xi,\tilde{\xi})\Big)=
\begin{cases}
\tilde{H}_{n_{\ell},(\xi_{0},\tilde{\xi}_{0})}(\xi,\tilde{\xi})  \;\;\;\;\;&\mbox{if $\ell^{'}=\ell$},\\
0                         \;\;\;\;\;&\mbox{if $\ell^{'}\neq \ell$}.
\end{cases}
\end{equation*}
Multiplying $\tilde{H}_{n_{\ell},(\xi_{0},\tilde{\xi}_{0})}(\xi,\tilde{\xi})$ on both sides of $(\ref{M12})$ gives
\begin{equation}\label{M13}
\tilde{H}_{n_{\ell},(\xi_{0},\tilde{\xi}_{0})}(\xi,\tilde{\xi})B_{2}(\xi,\tilde{\xi})^{\mathrm{T}}\hat{\Psi}_{2,(\xi_{0},\tilde{\xi}_{0})}(\xi,\tilde{\xi})=\ell\tilde{H}_{n_{\ell},(\xi_{0},\tilde{\xi}_{0})}(\xi,\tilde{\xi})D_{n_{\ell}}(\xi,\tilde{\xi})^{\mathrm{T}}\hat{\Psi}_{2,(\xi_{0},\tilde{\xi}_{0})}(\xi,\tilde{\xi}).
\end{equation}
Let $b_{2}$ be the Fourier coefficients of $B_{2}(\xi,\tilde{\xi})$, it follows from $(\ref{r3.16})$ and $\|\tilde{T}_{n_{\ell},(\xi_{0},\tilde{\xi}_{0})}\|\neq 0$ that
\begin{equation*}
\bigg\|\mathcal{F}^{-1}\Big(\ell\tilde{H}_{n_{\ell},(\xi_{0},\tilde{\xi}_{0})}(\xi,\tilde{\xi})D_{n_{\ell}}(\xi,\tilde{\xi})^{\mathrm{T}}\hat{\Psi}_{2,(\xi_{0},\tilde{\xi}_{0})}(\xi,\tilde{\xi})\Big)\bigg\|_{L^{p,q}}\geq \ell \|\tilde{T}_{n_{\ell},(\xi_{0},\tilde{\xi}_{0})}\|/2
\end{equation*}
and\\
\begin{equation*}
\bigg\|\mathcal{F}^{-1}\Big(\tilde{H}_{n_{\ell},(\xi_{0},\tilde{\xi}_{0})}(\xi,\tilde{\xi})B_{2}(\xi,\tilde{\xi})^{\mathrm{T}}\hat{\Psi}_{2,(\xi_{0},\tilde{\xi}_{0})}(\xi,\tilde{\xi})\Big)\bigg\|_{L^{p,q}}\leq  \|\tilde{T}_{n_{\ell},(\xi_{0},\tilde{\xi}_{0})}\|\;b_{2},
\end{equation*}
which contradicts with $(\ref{M13})$. This completes the proof of $(i)\Rightarrow (iii)$.
\subsection{Proof of $(iii)\Rightarrow (v)$}
Let $h_{(\eta_{\lambda_{1}},\eta_{\lambda_{2}})}(\xi,\tilde{\xi}),P_{(\eta_{\lambda_{1}},\eta_{\lambda_{2}})}(\xi,\tilde{\xi})$ and $\hat{\Psi}_{1,(\eta_{\lambda_{1}},\eta_{\lambda_{2}})}$ be as in Lemma $\ref{L8}$. Define
\begin{equation}\label{M14}
B_{(\eta_{\lambda_{1}},\eta_{\lambda_{2}})}(\xi,\tilde{\xi})=H_{(\eta_{\lambda_{1}},\eta_{\lambda_{2}})}(\xi,\tilde{\xi})\overline{P_{(\eta_{\lambda_{1}},\eta_{\lambda_{2}})}(\xi,\tilde{\xi})^{\mathrm{T}}}
\begin{pmatrix}
[\hat{\Psi}_{1,(\eta_{\lambda_{1}},\eta_{\lambda_{2}})},\hat{\Psi}_{1,(\eta_{\lambda_{1}},\eta_{\lambda_{2}})}](\xi,\tilde{\xi})^{-1}&0\\
0&I
\end{pmatrix}
P_{(\eta_{\lambda_{1}},\eta_{\lambda_{2}})}(\xi,\tilde{\xi}),
\end{equation}
where $H_{(\eta_{\lambda_{1}},\eta_{\lambda_{2}})}(\xi,\tilde{\xi})$ is a $2\pi$-periodic $C^{\infty}$ function such that $H_{(\eta_{\lambda_{1}},\eta_{\lambda_{2}})}(\xi,\tilde{\xi})=1$ on the support of $h_{(\eta_{\lambda_{1}},\eta_{\lambda_{2}})}$ and $H_{(\eta_{\lambda_{1}},\eta_{\lambda_{2}})}$ is supported in $B\Big((\eta_{\lambda_{1}},\eta_{\lambda_{2}}),2\delta_{(\eta_{\lambda_{1}},\eta_{\lambda_{2}})}\Big)+2\pi(\mathbf{Z}\times\mathbf{Z}^{d})$. Then the Fourier coefficients of $B_{(\eta_{\lambda_{1}},\eta_{\lambda_{2}})}(\xi,\tilde{\xi})$ belong to $\ell^{1,1}.$ Define $\Psi=(\psi_{1},\ldots,\psi_{r})^{\mathrm{T}}$ by
\begin{equation}\label{M15}
\hat{\Psi}(\xi,\tilde{\xi})=\sum\limits_{(\eta_{\lambda_{1}},\eta_{\lambda_{2}})\in \Lambda}h_{(\eta_{\lambda_{1}},\eta_{\lambda_{2}})}(\xi,\tilde{\xi})B_{(\eta_{\lambda_{1}},\eta_{\lambda_{2}})}(\xi,\tilde{\xi})\hat{\Phi}(\xi,\tilde{\xi}).
\end{equation}
Then $\Psi\in \mathcal{L}^{\infty,\infty}$ if $1< p,q<\infty$ and $\Psi\in \mathcal{W}(L^{1,1})$ if $p=1,\infty$ or $q=1,\infty$. For any $f\in V_{p,q}(\Phi)$, define
\begin{equation*}
g = \sum\limits_{i=1}^{r}\sum\limits_{j_{1}\in \mathbf{Z}}\sum\limits_{j_{2}\in \mathbf{Z}^{d}}\langle f,\psi_{i}(\cdot-j_{1},\cdot-j_{2})\rangle \phi_{i}(\cdot-j_{1},\cdot-j_{2}).
\end{equation*}
Then it sufficies to prove that $g=f$. By the definition of the space $V_{p,q}(\Phi),$ there exists $2\pi$-periodic distribution $D(\xi,\tilde{\xi})$ with Fourier coefficients belonging to $\ell^{p,q}$  such that
$\hat{f}(\xi,\tilde{\xi})=D(\xi,\tilde{\xi})^{\mathrm{T}}\hat{\Phi}(\xi,\tilde{\xi})$. Therefore, it follows from $(\ref{M14}),(\ref{M15})$ and Lemma $\ref{L8}$ that
\begin{eqnarray*}
\hat{g}(\xi,\tilde{\xi})=&&D(\xi,\tilde{\xi})^{\mathrm{T}}[\hat{\Phi},\hat{\Psi}](\xi,\tilde{\xi})\hat{\Phi}(\xi,\tilde{\xi})\\
=&&\sum\limits_{(\eta_{\lambda_{1}},\eta_{\lambda_{2}})\in \Lambda}h_{(\eta_{\lambda_{1}},\eta_{\lambda_{2}})}(\xi,\tilde{\xi})D(\xi,\tilde{\xi})^{\mathrm{T}}[\hat{\Phi},\hat{\Phi}](\xi,\tilde{\xi})\overline{B_{(\eta_{\lambda_{1}},\eta_{\lambda_{2}})}(\xi,\tilde{\xi})}\hat{\Phi}(\xi,\tilde{\xi})\\
=&&\sum\limits_{(\eta_{\lambda_{1}},\eta_{\lambda_{2}})\in \Lambda}h_{(\eta_{\lambda_{1}},\eta_{\lambda_{2}})}(\xi,\tilde{\xi})D(\xi,\tilde{\xi})^{\mathrm{T}}P_{(\eta_{\lambda_{1}},\eta_{\lambda_{2}})}(\xi,\tilde{\xi})^{-1}
\begin{pmatrix}
[\hat{\Psi}_{1,(\eta_{\lambda_{1}},\eta_{\lambda_{2}})},\hat{\Psi}_{1,(\eta_{\lambda_{1}},\eta_{\lambda_{2}})}](\xi,\tilde{\xi})&0\\
0&0\\
\end{pmatrix}\\
&&\times
\begin{pmatrix}
[\hat{\Psi}_{1,(\eta_{\lambda_{1}},\eta_{\lambda_{2}})},\hat{\Psi}_{1,(\eta_{\lambda_{1}},\eta_{\lambda_{2}})}](\xi,\tilde{\xi})^{-1}&0\\
0&I
\end{pmatrix}
\begin{pmatrix}
\hat{\Psi}_{1,(\eta_{\lambda_{1}},\eta_{\lambda_{2}})}(\xi,\tilde{\xi})\\
0
\end{pmatrix}\\
=&&\sum\limits_{(\eta_{\lambda_{1}},\eta_{\lambda_{2}})\in \Lambda}h_{(\eta_{\lambda_{1}},\eta_{\lambda_{2}})}(\xi,\tilde{\xi})D(\xi,\tilde{\xi})^{\mathrm{T}}P_{(\eta_{\lambda_{1}},\eta_{\lambda_{2}})}(\xi,\tilde{\xi})^{-1}
\begin{pmatrix}
\hat{\Psi}_{1,(\eta_{\lambda_{1}},\eta_{\lambda_{2}})}(\xi,\tilde{\xi})\\
0
\end{pmatrix}\\
=&&D(\xi,\tilde{\xi})^{\mathrm{T}}\hat{\Phi}(\xi,\tilde{\xi})=\hat{f}(\xi,\tilde{\xi}).
\end{eqnarray*}
Similarly, we can prove that $f= \sum\limits_{i=1}^{r}\sum\limits_{j_{1}\in \mathbf{Z}}\sum\limits_{j_{2}\in \mathbf{Z}^{d}}\langle f,\phi_{i}(\cdot-j_{1},\cdot-j_{2})\rangle \psi_{i}(\cdot-j_{1},\cdot-j_{2})$.
This completes the proof of $(iii)\Rightarrow (v)$.
\subsection{Proof of $(v)\Rightarrow (ii)$}
Let $f= \sum\limits_{i=1}^{r}\sum\limits_{j_{1}\in \mathbf{Z}}\sum\limits_{j_{2}\in \mathbf{Z}^{d}}\langle f,\phi_{i}(\cdot-j_{1},\cdot-j_{2})\rangle \psi_{i}(\cdot-j_{1},\cdot-j_{2}).$ Then it follows from Lemma $\ref{Lm 2.1}$ that\\
\begin{eqnarray}\label{r3.24}
\|f\|_{L^{p,q}}&&\leq\sum\limits_{i=1}^{r}\bigg\|\Big\{<f,\phi_{i}(\cdot-j_{1},\cdot-j_{2})>\Big\}_{(j_{1},j_{2})\in \mathbf{Z}\times\mathbf{Z}^{d}}\bigg\|_{\ell^{p,q}}\|\psi_{i}\|_{\mathcal{L}^{p,q}}\nonumber\\
&&\leq\bigg(\sum\limits_{i=1}^{r}\bigg\|\Big\{<f,\phi_{i}(\cdot-j_{1},\cdot-j_{2})>\Big\}_{(j_{1},j_{2})\in \mathbf{Z}\times\mathbf{Z}^{d}}\bigg\|_{\ell^{p,q}}\bigg)\bigg(\sum\limits_{i=1}^{r}\|\psi_{i}\|_{\mathcal{L}^{\infty,\infty}}\bigg).\;\;
\end{eqnarray}
By Lemma $\ref{L2.5}$, we can obtain
\begin{equation*}
\sum\limits_{i=1}^{r}\bigg\|\Big\{<f,\phi_{i}(\cdot-j_{1},\cdot-j_{2})>\Big\}_{(j_{1},j_{2})\in \mathbf{Z}\times\mathbf{Z}^{d}}\bigg\|_{\ell^{p,q}}\leq\bigg(\sum\limits_{i=1}^{r}\|\phi_{i}\|_{\mathcal{L}^{\infty,\infty}}\bigg)\|f\|_{L^{p,q}}.
\end{equation*}
Thus, the item $(ii)$ of Theorem $\ref{M3.1}$ holds for $A=\max\bigg\{\sum\limits_{i=1}^{r}\|\psi_{i}\|_{\mathcal{L}^{\infty,\infty}},\sum\limits_{i=1}^{r}\|\phi_{i}\|_{\mathcal{L}^{\infty,\infty}}\bigg\}$.
\subsection{Proof of $(ii)\Rightarrow (iii)$}
Let $k_{0}=\min\limits_{(\xi,\tilde{\xi})\in \mathbf{R}\times\mathbf{R}^{d}}rank\Big(\hat{\Phi}(\xi+2k_{1}\pi,\tilde{\xi}+2k_{2}\pi)\Big)_{(k_{1},k_{2})\in \mathbf{Z}\times\mathbf{Z}^{d}}$ and $\Omega_{k_{0}}=\Big\{(\xi,\tilde{\xi}):rank\Big(\hat{\Phi}(\xi+2k_{1}\pi,\tilde{\xi}+2k_{2}\pi)\Big)_{(k_{1},k_{2})\in \mathbf{Z}\times\mathbf{Z}^{d}}>k_{0}\Big\}$. By Lemma $\ref{L7}$, it suffices to prove that $\Omega_{k_{0}}=\varnothing$. On the contrary, suppose that $\Omega_{k_{0}}\neq \varnothing$. Let $(\xi_{0},\tilde{\xi}_{0})\in \partial\Omega_{k_{0}},\:\Psi_{1,(\xi_{0},\tilde{\xi}_{0})},\;\Psi_{2,(\xi_{0},\tilde{\xi}_{0})},\;P_{(\xi_{0},\tilde{\xi}_{0})}(\xi,\tilde{\xi}),\;H_{n,(\xi_{0},\tilde{\xi}_{0})}(\xi,\tilde{\xi}),$\\
$\tilde{H}_{n,(\xi_{0},\tilde{\xi}_{0})}(\xi,\tilde{\xi})$ and $\delta_{(\xi_{0},\tilde{\xi}_{0})}$ be as in the proof of $(i)\Rightarrow (iii)$. It follows from $(\ref{M2})$ and the continuity of $[\hat{\Psi}_{1,(\xi_{0},\tilde{\xi}_{0})},\hat{\Psi}_{1,(\xi_{0},\tilde{\xi}_{0})}](\xi,\tilde{\xi})$ that we can chose a $n_{0}$ such that $2^{-n_{0}}<\delta_{(\xi_{0},\tilde{\xi}_{0})}$ and
\begin{eqnarray*}
\alpha_{n}(\xi,\tilde{\xi})=&&[\hat{\Psi}_{1,(\xi_{0},\tilde{\xi}_{0})},\hat{\Psi}_{1,(\xi_{0},\tilde{\xi}_{0})}](\xi_{0},\tilde{\xi}_{0})\\
&&+H_{n,(\xi_{0},\tilde{\xi}_{0})}(\xi,\tilde{\xi})\bigg([\hat{\Psi}_{1,(\xi_{0},\tilde{\xi}_{0})},\hat{\Psi}_{1,(\xi_{0},\tilde{\xi}_{0})}](\xi,\tilde{\xi})-[\hat{\Psi}_{1,(\xi_{0},\tilde{\xi}_{0})},\hat{\Psi}_{1,(\xi_{0},\tilde{\xi}_{0})}](\xi_{0},\tilde{\xi}_{0})\bigg)
\end{eqnarray*}
is nonsingular for all $n\geq n_{0}$. Given any $(r-k_{0})\times 1$ $2\pi$-periodic distribution matrix $F(\xi,\tilde{\xi})$ with Fourier coefficients belonging to $\ell^{p,q}$, define $g_{n},n\geq n_{0}+1$, by
\begin{eqnarray}\label{r22}
\hat{g}_{n}(\xi,\tilde{\xi})=&&\tilde{H}_{n,(\xi_{0},\tilde{\xi}_{0})}(\xi,\tilde{\xi})\bigg(-F(\xi,\tilde{\xi})^{\mathrm{T}}[\hat{\Psi}_{2,(\xi_{0},\tilde{\xi}_{0})},\hat{\Psi}_{1,(\xi_{0},\tilde{\xi}_{0})}](\xi,\tilde{\xi})\Big(\alpha_{n}(\xi,\tilde{\xi})\Big)^{-1},F(\xi,\tilde{\xi})^{\mathrm{T}}\bigg)\nonumber\\
&&\times
\begin{pmatrix}
\hat{\Psi}_{1,(\xi_{0},\tilde{\xi}_{0})}(\xi,\tilde{\xi})\\
\hat{\Psi}_{2,(\xi_{0},\tilde{\xi}_{0})}(\xi,\tilde{\xi})
\end{pmatrix}.\;\;\;\;\;\;\;\;\;
\end{eqnarray}
Then $g_{n}\in V_{p,q}(\Phi)$ and
\begin{eqnarray}\label{M16}
[\hat{g}_{n},\hat{\Psi}_{1,(\xi_{0},\tilde{\xi}_{0})}](\xi,\tilde{\xi})=&&\tilde{H}_{n,(\xi_{0},\tilde{\xi}_{0})}(\xi,\tilde{\xi})\bigg(-F(\xi,\tilde{\xi})^{\mathrm{T}}[\hat{\Psi}_{2,(\xi_{0},\tilde{\xi}_{0})},\hat{\Psi}_{1,(\xi_{0},\tilde{\xi}_{0})}](\xi,\tilde{\xi})\Big(\alpha_{n}(\xi,\tilde{\xi})\Big)^{-1},F(\xi,\tilde{\xi})^{\mathrm{T}}\bigg)\nonumber\\
&&\times
\begin{pmatrix}
[\hat{\Psi}_{1,(\xi_{0},\tilde{\xi}_{0})},\hat{\Psi}_{1,(\xi_{0},\tilde{\xi}_{0})}](\xi,\tilde{\xi})\\
[\hat{\Psi}_{2,(\xi_{0},\tilde{\xi}_{0})},\hat{\Psi}_{1,(\xi_{0},\tilde{\xi}_{0})}](\xi,\tilde{\xi})
\end{pmatrix}\nonumber\\
=&&0,
\end{eqnarray}
where we have used $(\ref{M8})$ and the fact that
$\alpha_{n}(\xi,\tilde{\xi})=[\hat{\Psi}_{1,(\xi_{0},\tilde{\xi}_{0})},\hat{\Psi}_{1,(\xi_{0},\tilde{\xi}_{0})}](\xi,\tilde{\xi})$
on the support of $\tilde{H}_{n,(\xi_{0},\tilde{\xi}_{0})}.$
Moreover, by Lemma $\ref{L2.5},(\ref{M8}),(\ref{M16})$ and
\begin{equation*}
P_{(\xi_{0},\tilde{\xi}_{0})}(\xi,\tilde{\xi})\hat{\Phi}(\xi,\tilde{\xi})=
\begin{pmatrix}
\hat{\Psi}_{1,(\xi_{0},\tilde{\xi}_{0})}(\xi,\tilde{\xi})\\
\hat{\Psi}_{2,(\xi_{0},\tilde{\xi}_{0})}(\xi,\tilde{\xi})
\end{pmatrix},
\end{equation*}
we obtain
\begin{eqnarray}\label{M17}
\ \ \ \ \|[\hat{g}_{n},\hat{\Phi}](\xi,\tilde{\xi})\|_{\ell_{\ast}^{p,q}}
&&=\bigg\|\bigg[\hat{g}_{n},
\begin{pmatrix}
\hat{\Psi}_{1,(\xi_{0},\tilde{\xi}_{0})}\nonumber\\
\hat{\Psi}_{2,(\xi_{0},\tilde{\xi}_{0})}
\end{pmatrix}\bigg]
(\xi,\tilde{\xi})\overline{P_{(\xi_{0},\tilde{\xi}_{0})}(\xi,\tilde{\xi})^{-\mathrm{T}}}\bigg\|_{\ell_{\ast}^{p,q}}\nonumber\\
&&\leq\bigg\|\bigg[\hat{g}_{n},
\begin{pmatrix}
\hat{\Psi}_{1,(\xi_{0},\tilde{\xi}_{0})}\nonumber\\
\hat{\Psi}_{2,(\xi_{0},\tilde{\xi}_{0})}
\end{pmatrix}\bigg]
(\xi,\tilde{\xi})\bigg\|_{\ell_{\ast}^{p,q}}\|\overline{P_{(\xi_{0},\tilde{\xi}_{0})}(\xi,\tilde{\xi})^{-\mathrm{T}}}\|_{\ell_{\ast}^{1,1}}\nonumber\\
&&\leq C\bigg\|\bigg[\hat{g}_{n},
\begin{pmatrix}
\hat{\Psi}_{1,(\xi_{0},\tilde{\xi}_{0})}\nonumber\\
\hat{\Psi}_{2,(\xi_{0},\tilde{\xi}_{0})}
\end{pmatrix}\bigg]
(\xi,\tilde{\xi})\bigg\|_{\ell_{\ast}^{p,q}}\nonumber\\
&&= C\Big\|[\hat{g}_{n},\hat{\Psi}_{2,(\xi_{0},\tilde{\xi}_{0})}](\xi,\tilde{\xi})\Big\|_{\ell_{\ast}^{p,q}}\;\;\;\;\;\;\;\;\;\;\;\nonumber\\
&&= C\Big\|[\hat{g}_{n},H_{n,(\xi_{0},\tilde{\xi}_{0})}\hat{\Psi}_{2,(\xi_{0},\tilde{\xi}_{0})}](\xi,\tilde{\xi})\Big\|_{\ell_{\ast}^{p,q}}\;\;\;\;\;\;\;\;\;\;\;\nonumber\\
&&\leq  C\|g_{n}\|_{L^{p,q}}\Big\|\mathcal{F}^{-1}(H_{n,(\xi_{0},\tilde{\xi}_{0})}\hat{\Psi}_{2,(\xi_{0},\tilde{\xi}_{0})})\Big\|_{\mathcal{L}^{\infty,\infty}},
\end{eqnarray}
where $\|\overline{P_{(\xi_{0},\tilde{\xi}_{0})}(\xi,\tilde{\xi})^{-\mathrm{T}}}\|_{\ell_{\ast}^{1,1}}\leq C$.
It follows from Lemma $\ref{L9}$ that
\begin{eqnarray*}
&&\Big\|\mathcal{F}^{-1}(H_{n,(\xi_{0},\tilde{\xi}_{0})}\hat{\Psi}_{2,(\xi_{0},\tilde{\xi}_{0})})\Big\|_{\mathcal{L}^{\infty,\infty}}\\
&&=\Big\|\sum\limits_{j_{1}\in\mathbf{Z}}\sum\limits_{j_{2}\in\mathbf{Z}^{d}}h_{n,(\xi_{0},\tilde{\xi}_{0})}(j_{1},j_{2})\Psi_{2,(\xi_{0},\tilde{\xi}_{0})}(\cdot-j_{1},\cdot-j_{2})\Big\|_{\mathcal{L}^{\infty,\infty}}\\
&&\leq2^{-n(d+1)}\Big\|\sum\limits_{j_{1}\in\mathbf{Z}}\sum\limits_{j_{2}\in\mathbf{Z}^{d}} \hat{H}(2^{-n}j_{1},2^{-n}j_{2})e^{-i(\xi_{0},\tilde{\xi}_{0})\cdot(\cdot+j_{1},\cdot+j_{2})}\Psi_{2,(\xi_{0},\tilde{\xi}_{0})}(\cdot+j_{1},\cdot+j_{2})\Big\|_{\mathcal{L}^{\infty,\infty}}.
\end{eqnarray*}
Then
\begin{equation*}
\lim\limits_{n\rightarrow\infty}\Big\|\mathcal{F}^{-1}(H_{n,(\xi_{0},\tilde{\xi}_{0})}\hat{\Psi}_{2,(\xi_{0},\tilde{\xi}_{0})})\Big\|_{\mathcal{L}^{\infty,\infty}}=0.
\end{equation*}
This together with $(\ref{M17})$ leads to the existence of $\epsilon_{n},n\geq n_{0}$, such that
\begin{equation}\label{M18}
\|[\hat{g}_{n},\hat{\Phi}](\xi,\tilde{\xi})\|_{\ell_{\ast}^{p,q}}\leq \epsilon_{n}\|g_{n}\|_{L^{p,q}}
\end{equation}
and $\lim\limits_{n\rightarrow \infty}\epsilon_{n}=0$. On the other hand, by the assumption $(ii)$ we obtain
\begin{equation*}
\|[\hat{g}_{n},\hat{\Phi}](\xi,\tilde{\xi})\|_{\ell_{\ast}^{p,q}}=\bigg\|\bigg\{\int_{\mathbf{R}}\int_{\mathbf{R}^{d}}g_{n}(x_{1},x_{2})\overline{\Phi(x_{1}-j_{1},x_{2}-j_{2})}\bigg\}_{(j_{1},j_{2})\in\mathbf{Z}\times\mathbf{Z}^{d}}\bigg\|_{\ell^{p,q}}\geq A^{-1}\|g_{n}\|_{L^{p,q}}.
\end{equation*}
This together with $(\ref{M18})$ proves that there exists an integer $n_{1}\geq n_{0}+1$ such that
\begin{equation*}
g_{n}\equiv 0\;\;\;for\; n\geq n_{1}.
\end{equation*}
Thus for any $(r-k_{0})\times 1$ $2\pi$-periodic distribution matrix $F(\xi,\tilde{\xi})$ with Fourier coefficients belonging to $\ell^{p,q}$, we can obtain
\begin{eqnarray*}
&&\tilde{H}_{n,(\xi_{0},\tilde{\xi}_{0})}(\xi,\tilde{\xi})F(\xi,\tilde{\xi})^{\mathrm{T}}[\hat{\Psi}_{2,(\xi_{0},\tilde{\xi}_{0})},\hat{\Psi}_{1,(\xi_{0},\tilde{\xi}_{0})}](\xi,\tilde{\xi})\Big(\alpha_{n}(\xi,\tilde{\xi})\Big)^{-1}\hat{\Psi}_{1,(\xi_{0},\tilde{\xi}_{0})}(\xi,\tilde{\xi})\\
&&=\tilde{H}_{n,(\xi_{0},\tilde{\xi}_{0})}(\xi,\tilde{\xi})F(\xi,\tilde{\xi})^{\mathrm{T}}\hat{\Psi}_{2,(\xi_{0},\tilde{\xi}_{0})}(\xi,\tilde{\xi})
\end{eqnarray*}
from $(\ref{r22})$ for any $n\geq n_{1}$. Hence
\begin{eqnarray}\label{M20}
&&\tilde{H}_{n,(\xi_{0},\tilde{\xi}_{0})}(\xi,\tilde{\xi})[\hat{\Psi}_{2,(\xi_{0},\tilde{\xi}_{0})},\hat{\Psi}_{1,(\xi_{0},\tilde{\xi}_{0})}](\xi,\tilde{\xi})\Big(\alpha_{n}(\xi,\tilde{\xi})\Big)^{-1}\hat{\Psi}_{1,(\xi_{0},\tilde{\xi}_{0})}(\xi,\tilde{\xi})\nonumber\\
&&=\tilde{H}_{n,(\xi_{0},\tilde{\xi}_{0})}(\xi,\tilde{\xi})\hat{\Psi}_{2,(\xi_{0},\tilde{\xi}_{0})}(\xi,\tilde{\xi}).
\end{eqnarray}
This together with $(\ref{M2})$ and $(\ref{M4})$ leads to\\
\begin{equation*}
\tilde{H}_{n,(\xi_{0},\tilde{\xi}_{0})}(\xi,\tilde{\xi})[\hat{\Psi}_{2,(\xi_{0},\tilde{\xi}_{0})},\hat{\Psi}_{1,(\xi_{0},\tilde{\xi}_{0})}](\xi,\tilde{\xi})\Big(\alpha_{n}(\xi,\tilde{\xi})\Big)^{-1}=0,\;\;\;(\xi,\tilde{\xi})\in B\Big((\xi_{0},\tilde{\xi}_{0}),2^{-n_{1}}\Big)+2\pi(\mathbf{Z}\times\mathbf{Z}^{d}).
\end{equation*}
Substituting this into $(\ref{M20})$, we obtain
\begin{equation*}
\hat{\Psi}_{2,(\xi_{0},\tilde{\xi}_{0})}(\xi,\tilde{\xi})\equiv 0\;\;\;on\;B\Big((\xi_{0},\tilde{\xi}_{0}),2^{-n_{1}}\Big)+2\pi(\mathbf{Z}\times\mathbf{Z}^{d}),
\end{equation*}
which contradicts with $(\ref{M5})$. This completes the proof of $(ii)\Rightarrow (iii)$.\\

\noindent \textbf{Acknowledgement }\ \ The project is partially  supported by the National Natural Science Foundation of China (Nos. 11661024, 11671107) and the Guangxi Natural Science Foundation (Nos. 2019GXNSFFA245012, 2017GXNSFAA198194), Guangxi Key Laboratory of Cryptography and Information Security (No. GCIS201614), Guangxi Colleges and Universities Key Laboratory of Data Analysis and Computation.

 \end{document}